\documentclass[final,1p,times]{elsarticle3}

 \usepackage{amssymb}
 \usepackage{amsthm}

\usepackage{amsmath}
\usepackage{amsfonts}
\usepackage{latexsym}
\usepackage{mathrsfs}
\usepackage{times}
\usepackage{array}
\usepackage{amscd}
\usepackage{a4}
\usepackage{pstricks}

\input{boxedeps} 
\SetEPSFDirectory{} 
\HideDisplacementBoxes
\SetRokickiEPSFSpecial  
\newtheorem{theorem}{Theorem}[section]

\newdefinition{example}{Example}[section]
\newdefinition{remark}{Remark}[section]
\newdefinition{aside}{Aside}[section]
\input xy 
\xyoption{all} 
\xyoption{2cell} 
\DeclareMathAlphabet{\ams}{U}{msb}{m}{n}
\DeclareMathAlphabet{\goth}{U}{euf}{m}{n}

\def\sl{S\kern-1pt L}

\def\ml{M\kern-1pt L}

\def\sp{S\kern-1.5pt p}

\def\ov{\overline}

\def\HH{\mathcal H}

\def\gS{\goth{S}}

\def\ve{\varepsilon}

\def\aa{\alpha}

\def\bb{\beta}
\def\ss{\sigma}

\def\ve{\varepsilon}

\def\Z{\ams{Z}}
\def\R{\ams{R}}
\def\Q{\ams{Q}}
\def\F{\ams{F}}

\def\quo{/\kern -.45em\sim}

\def\blob{\bullet}

\def\Langle{\langle\kern -2pt\langle}
\def\Rangle{\rangle\kern -1.9pt\rangle}

\newpsobject{showgrid}{psgrid}{subgriddiv=1,griddots=10,gridlabels=6pt,gridcolor=red}
\journal{Expositiones Mathematicae}

\begin{document}

\begin{frontmatter}



\title{A (very short) introduction to buildings\tnoteref{label1}}
\tnotetext[label1]{Based
on a short series of lectures given at the
Universit\'e de Fribourg, Switzerland in June 2011. The author is grateful for the
department's hospitality both then and over the years, and
particularly to Ruth Kellerhals for her many kindnesses. 
He also thanks Laura Chiobanu, Paul Turner, and was partially
supported by 
Swiss National Science Foundation
grant 200021-131967
and Marie Curie Reintegration Grant 230889. He also thanks the referees
for a number of helpful suggestions and references.}


\author{Brent Everitt}
\ead{brent.everitt@york.ac.uk}
\address{Department of Mathematics, University of York, York
YO10 5DD, United Kingdom.}

\begin{abstract}
This is an informal elementary
introduction to buildings -- what they are and where they come
from.
\end{abstract}




\end{frontmatter}




This is an informal elementary
introduction to buildings -- written for, and by,  a non-expert. 
The aim is to
get to the definition of a building and feel that it is an entirely
natural thing. To maintain the lecture style examples have replaced
proofs. The notes at the end indicate where these proofs can
be found. 

Most of what we say has its origins in the work of Jacques Tits, and
our account borrows heavily from
the books of Abramenko and Brown \cite{Abramenko_Brown08} and
Ronan \cite{Ronan09}. Section \ref{lecture1} illustrates all the
essential features of a
building in the context of an example, but without mentioning any
building terminology. In principle anyone could read
this. Sections \ref{lecture2}-\ref{lecture4} firm-up and generalize these
specifics: Coxeter groups appear in \S \ref{lecture2}, chambers systems in
\S \ref{lecture3} and the definition of a building in \S
\ref{lecture4}. Section \ref{lecture5} 
addresses where buildings come from by describing
the first important example: the spherical building of an algebraic group.

\section{The flag complex of a vector space}
\label{lecture1}

\begin{figure}
  \centering
\begin{pspicture}(0,0)(12,6)
\rput(3,3){
\rput(0,0){\BoxedEPSF{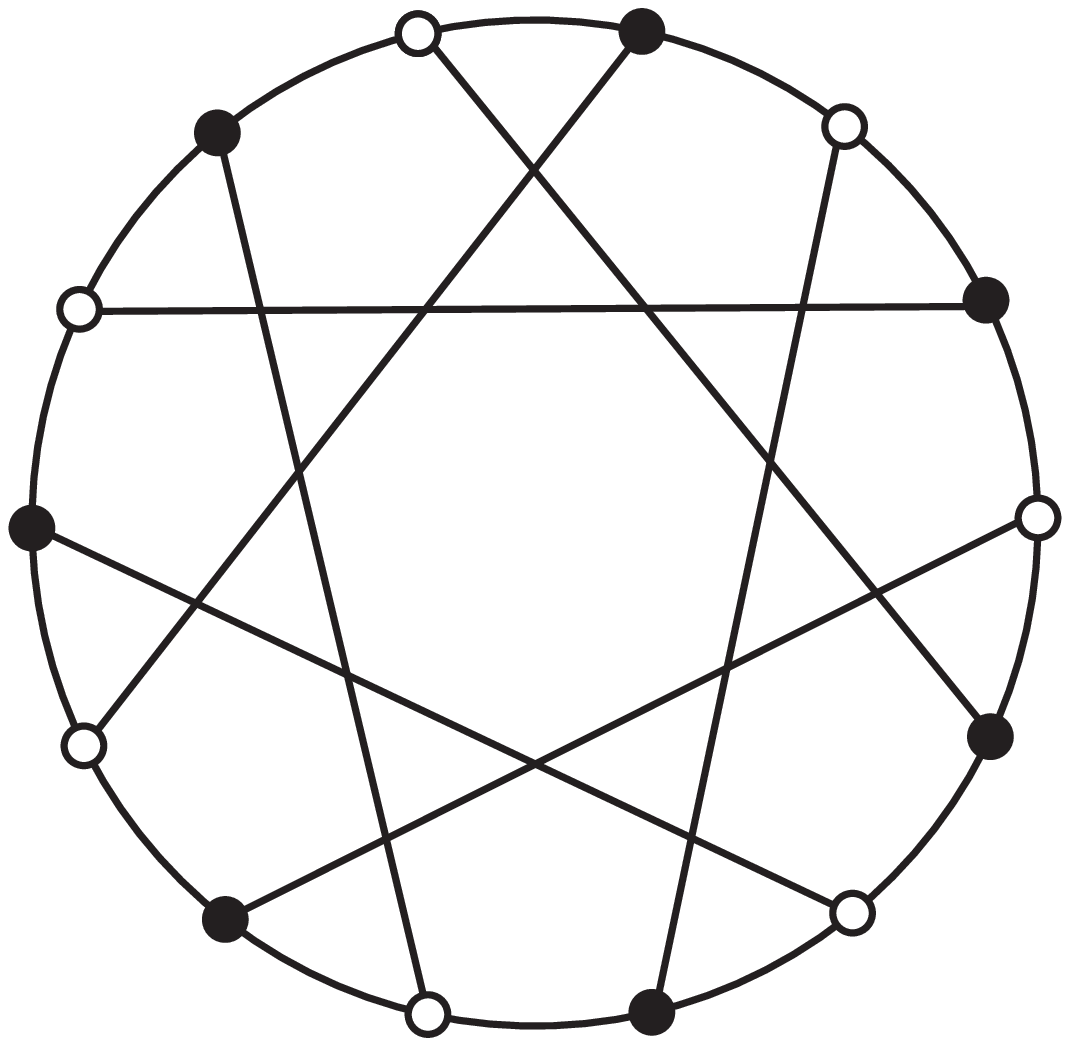 scaled 500}}
}
\rput(9,3){
\rput(0,0){\BoxedEPSF{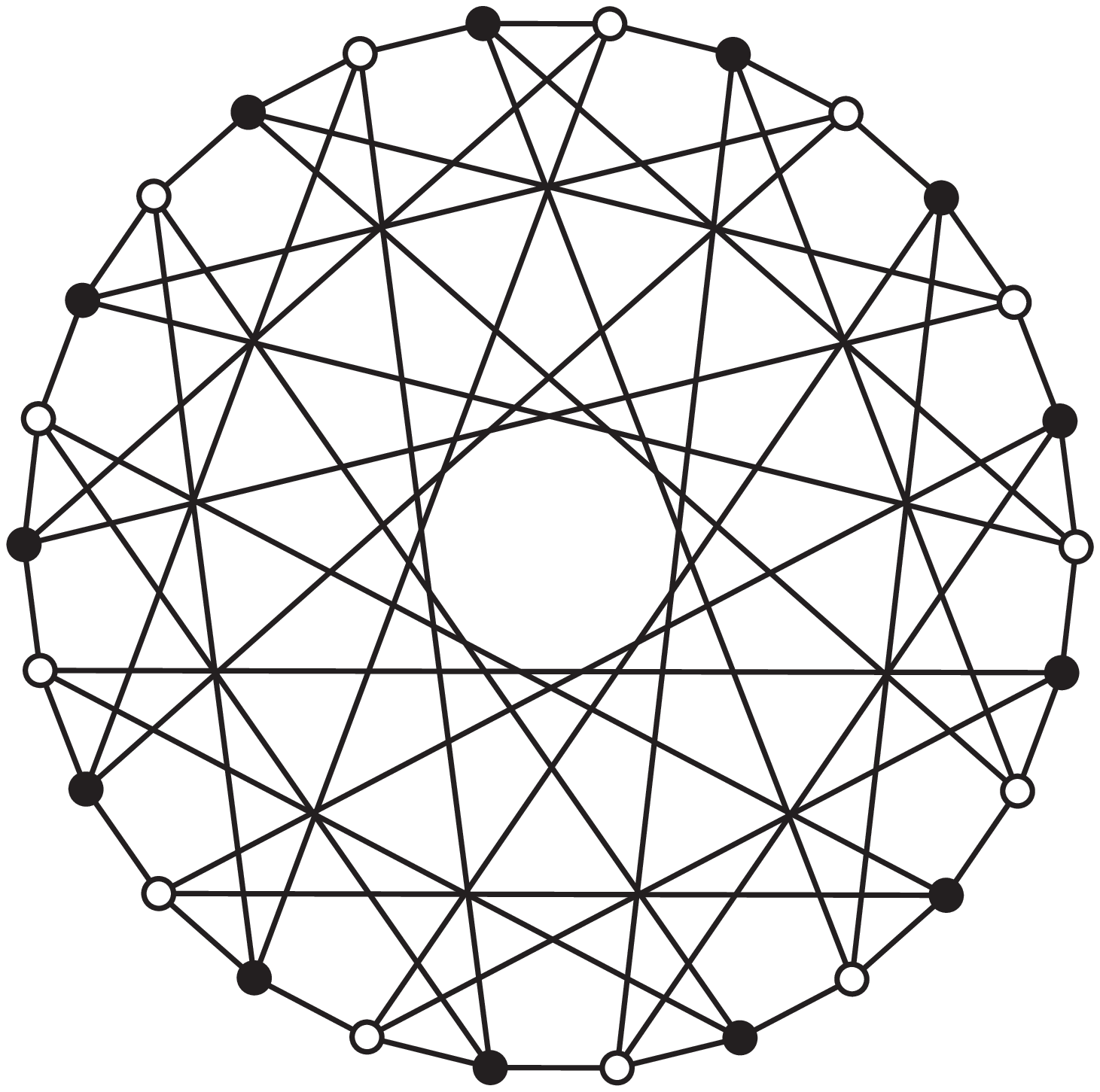 scaled 375}}
}
\end{pspicture}
\caption{The flag complex $\Delta$ of the three dimensional vector space
  over the fields of order $2$ \emph{(left)} and $3$ \emph{(right)}.}
  \label{fig:flagcomplex}
\end{figure}

Let $V$ be a three dimensional vector space over a field $k$. Let
$\Delta$ be the graph with vertices 
the non-trivial proper subspaces of $V$, and an edge connecting
the vertices $V_i$ and
$V_j$ whenever $V_i$ is a subspace of $V_j$\/:
$$
\xymatrix{
\begin{pspicture}(0,0)(2.65,0.1)
\rput(0,-0.4){
\psline[](0.5,.5)(2,.5)
\pscircle[fillstyle=solid,fillcolor=black](2,.5){0.15}
\pscircle[fillstyle=solid,fillcolor=white](0.5,.5){0.15}
\rput(0.05,.5){$V_i$}\rput(2.45,.5){$V_j$}}
\end{pspicture}
\ar@{<=>}[r]
&  V_i\subset V_j.
}
$$
Figure \ref{fig:flagcomplex} shows the graph $\Delta$ when $k$ is
the field of orders $q=2$ and $3$. There are $1+q+q^2$ one dimensional subspaces --
illustrated by the white vertices -- and $1+q+q^2$ two dimensional
subspaces, illustrated by the black vertices. Each one dimensional
space is contained in $1+q$ two dimensional spaces and each
two dimensional space contains $1+q$ one dimensional spaces. The
duality here might remind the reader of projective
geometry.
Call the edges $V_i\subset V_j$ of $\Delta$ \emph{chambers\/}.

Some more structure can be wrung out of this picture:
there is an ``$\gS_3$-valued metric'', with $\gS_3$ the symmetric
group, 
that gives the shortest route(s)
through $\Delta$ between any two chambers.
To see how, suppose $c,c'$ are chambers and we want
a shortest route of edges connecting them:
$$
\xymatrix{
c=V_1\subset V_2  \ar@{~>}[rr]^-{\rm{shortest\,\,route}}
&&  c'=V_1'\subset V_2'.
}
$$
Make $c$ and $c'$ as different as
possible by assuming that $V_1\not=V_1'$, $V_2\not=V_2'$ and
$V_2\cap V_2'$ is a line different from $V_1,V_1'$. Changing
notation, let $L_1,L_2,L_3$ be lines
with $L_1=V_1$, $L_3=V_1'$ and $L_2=V_2\cap V_2'$.
One then gets $V_2=L_1+L_2$
and $V_2'=L_2+L_3$. 

\begin{figure}
  \centering
\begin{pspicture}(0,0)(12,5.8)
\rput(1.25,1.5){
\rput(1.5,1.5){\BoxedEPSF{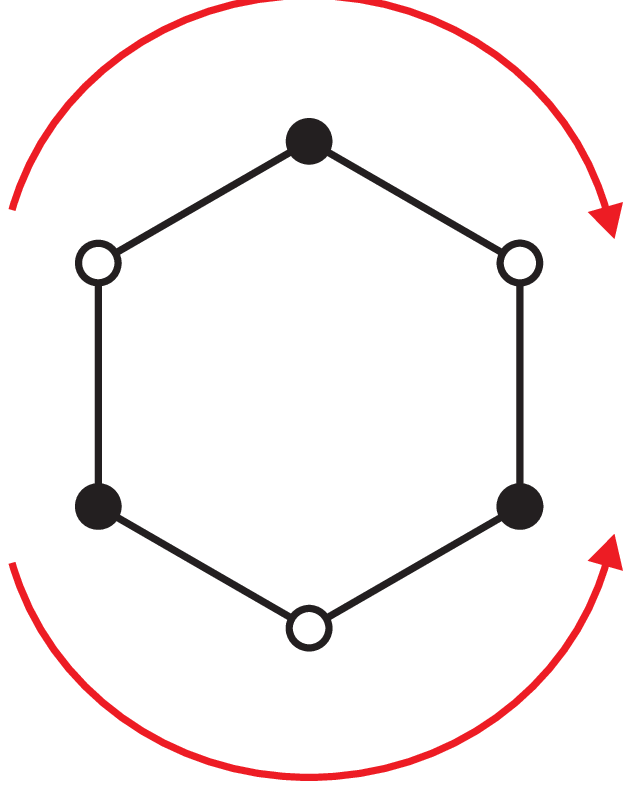 scaled 500}}
\rput(1.5,3){$\scriptstyle{L_1+L_3}$}
\rput(3.15,0.9){$\scriptstyle{L_2+L_3}$}
\rput(-0.15,0.9){$\scriptstyle{L_1+L_2}$}
\rput(0.1,2.1){$\scriptstyle{L_1}$}
\rput(1.5,-0.05){$\scriptstyle{L_2}$}
\rput(2.9,2.1){$\scriptstyle{L_3}$}
\rput(0.2,1.5){${\red c}$}\rput(2.8,1.5){${\red c'}$}
\rput(1.5,2.45){$\scriptstyle{1}$}
\rput(2.3,2.05){$\scriptstyle{2}$}
\rput(2.3,0.95){$\scriptstyle{1}$}
\rput(1.5,0.55){$\scriptstyle{2}$}
\rput(0.7,0.95){$\scriptstyle{1}$}
\rput(0.7,2.05){$\scriptstyle{2}$}
\rput(1.5,3.7){$\scriptstyle{s_2s_1s_2}$}
\rput(1.5,-0.75){$\scriptstyle{s_1s_2s_1}$}
}
\rput(-0.75,0){
\rput(9.5,3){
\rput(0,0){\BoxedEPSF{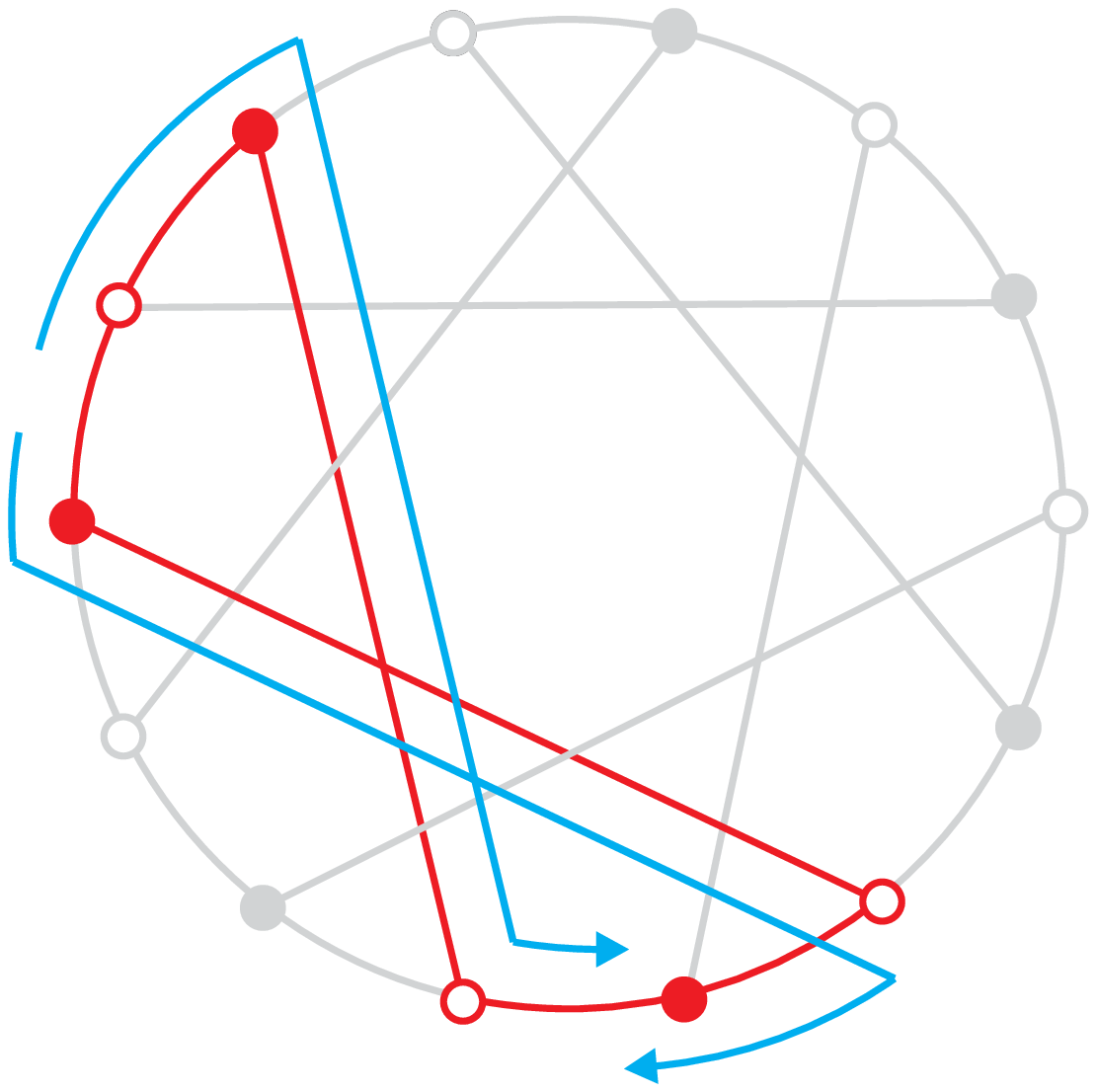 scaled 500}}
}
\rput(-0.5,0){\rput(7.3,3.8){${\red c}$}\rput(10.1,0.4){${\red c'}$}}
}
\end{pspicture}
\caption{Local picture of $\Delta$ containing the pair of chambers
$c,c'$ and the shortest routes between them \emph{(left)\/}; situating the pair $c,c'$ in a local picture
with the shortest routes
$s_1s_2s_1=s_2s_1s_2=(1,3)$ \emph{(right)\/}.}
  \label{fig:localpicture}
\end{figure}

We get a small piece of $\Delta$,
a local picture containing $c,c'$, as in Figure \ref{fig:localpicture}. The
field $k$ wasn't mentioned at all in the previous paragraph, so this is the
local picture for $\Delta$ over any field. The
\emph{global\/} picture gets more complicated however as the field $k$
gets bigger as Figure \ref{fig:flagcomplex} illustrates.

Say that chambers are \emph{$i$-adjacent\/} if any
difference between them occurs only in the $i$-th position, so $V_1\subset
V_2\supset V_1',\, (V_1\not=V_1')$ are a pair of $1$-adjacent chambers and $V_2\supset
V_1\subset V_2',\,(V_2\not= V_2')$ a pair of $2$-adjacent chambers (a
chamber is also $i$-adjacent to itself for any $i$). Place the label $i$
on a vertex of the local picture in Figure \ref{fig:localpicture} if the two chambers meeting at the
vertex are $i$-adjacent.

The shortest routes from $c$ to $c'$ \emph{in the local picture\/} are
given by
$$
\xymatrix{
c\ar[r]^-{s_2s_1s_2}_-{s_1s_2s_1}
&  c'
}
$$
where the route $s_1s_2s_1$ means cross a $1$-labeled vertex, then a
$2$-labeled vertex and then a $1$-labeled vertex. Routes are read 
from left to right, although it obviously doesn't matter with
the two above. These routes then take values in the symmetric group $\gS_3$ by letting $s_1=(1,2)$ and
$s_2=(2,3)$, so that both $s_1s_2s_1$ and $s_2s_1s_2$ give the
permutation $(1,3)\in\gS_3$. 
Our actions will always be on the left, so in particular permutations
in $\gS_3$ are composed from right to left.
Define the $\gS_3$-distance between
$c,c'$ to be $\delta(c,c')=(1,3)$.

For an arbitrary pair of chambers define $\delta(c,c')$
to be the element of $\gS_3$ obtained by situating the chambers $c,c'$
in some local picture and taking the shortest
route(s) as in Figure \ref{fig:localpicture}. The resulting map
$\delta:\Delta\times\Delta\rightarrow
\gS_3$ can be thought of as a metric on $\Delta$ taking values in $\gS_3$.

We will see in \S \ref{lecture4}
why this map is well defined and doesn't depend on which local
picture we choose containing $c,c'$,
although an ad-hoc argument
shows that an element of $\gS_3$ can be associated
in a canonical fashion to
any pair of chambers. Take the $c,c'$ above and write
$$
c=0\subset L_1\subset L_1+L_2\subset V=V_0\subset V_1\subset
V_2\subset V_3
$$
and $c'=V'_0\subset\cdots\subset V'_3$ similarly. For each $i$ 
the filtration $V_0\subset V_1\subset V_2\subset V_3$ of $V$
induces a
filtration 
of the one dimensional quotient $V'_i/V'_{i-1}$:
\begin{equation}
  \label{eq:13}
(V'_i\cap V_0)/V'_{i-1}\subset\cdots\subset (V'_i\cap V_3)/V'_{i-1}  
\end{equation}
(by $(V'_i\cap V_0)/V'_{i-1}$, etc, we mean the image of
  $V'_i\cap V_0$ under the quotient map $V\rightarrow V/V'_{i-1}$).
Any filtration of a one dimensional space
must start with a sequence of trivial subspaces and end with
a sequence of $V'_i/V'_{i-1}$'s. At some point in the middle the
filtration jumps from being zero dimensional to one dimensional; for
the $c,c'$ above:
$$
\begin{tabular}{cccc}
 $i$&$V'_i/V'_{i-1}$&filtration (\ref{eq:13})&``jump index" $j$\\
\hline
$1$&$L_3$&$0\subset 0\subset 0\subset L_3$&$3$\\
$2$&$(L_2+L_3)/L_3$&$0\subset 0\subset (L_2+L_3)/L_3\subset (L_2+L_3)/L_3$&$2$\\
$3$&$V/(L_2+L_3)$&$0\subset V/(L_2+L_3)\subset V/(L_2+L_3)\subset V/(L_2+L_3)$&$1$\\
\hline
\end{tabular}
$$
Defining $\pi(i)=j$ gives $\pi=(1,3)\in\gS_3$. 
Summarizing:

\paragraph{First rough definition of a building} A building is a set of
\emph{chambers\/} with $i$-adjacency between them, the $i$ coming from some set
$S$, together with a ``$W$-valued metric'' for $W$ some group.

\paragraph{}Returning to the running example, the symmetric group $\gS_3$ is
a reflection group, with Figure \ref{fig:localpicture}
and the resulting metric $\delta$ coming from the geometry of these
reflections. 
To see why suppose we have a three dimensional Euclidean space -- a real vector space
with an inner product. Let $v_1,v_2,v_3$ be
an orthonormal basis and let $\gS_3$ act on the space
by permuting coordinates: $\pi\cdot v_i:=v_{\pi(i)}$ for 
$\pi\in\gS_3$ (and
extend linearly). This action is not
essential as the vector $v=v_1+v_2+v_3$ is fixed by all
$\pi\in\gS_3$. This can be gotten around by passing to the perp
space 
$$
{\textstyle v^\perp=\{\sum\lambda_i v_i\,|\,\sum\lambda_i=0\}}.
$$
The picture to keep in mind is the following,
where $v^\perp$ is translated off the origin to make it easier to
see:
$$
\begin{pspicture}(0,0)(14.3,3.5)
\rput(0,0){
\rput(4.6,1.75){
\rput(0,0){\BoxedEPSF{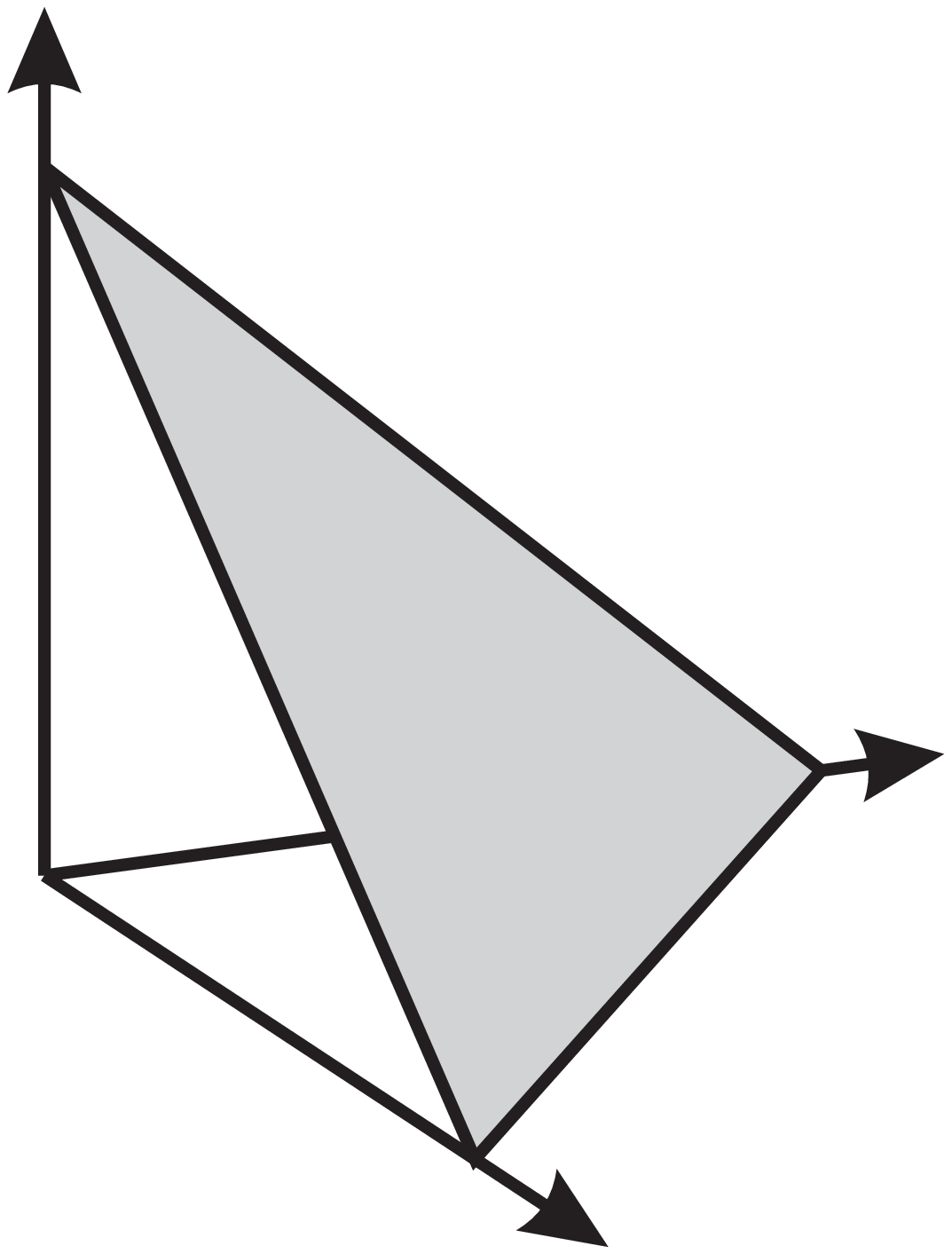 scaled 225}}
}
\rput(4.1,0.4){$v_1$}\rput(5.75,1.7){$v_2$}\rput(3.1,2.5){$v_3$}
\rput(5.7,0.6){${\red s_1}$}\rput(4.1,3.3){${\red s_1}$}
\rput(5,2.5){$v^\perp+\frac{1}{3}v$}
\psbezier[showpoints=false,linecolor=red]{->}(3.3,3.2)(3,3.75)(4,3.75)(3.7,3.2)
\psbezier[showpoints=false,linecolor=red]{<->}(5,0.4)(5.5,0.5)(5.6,1)(5.5,1.2)
}
\rput(0,0){
\rput(9.2,1.75){
\rput(0,0){\BoxedEPSF{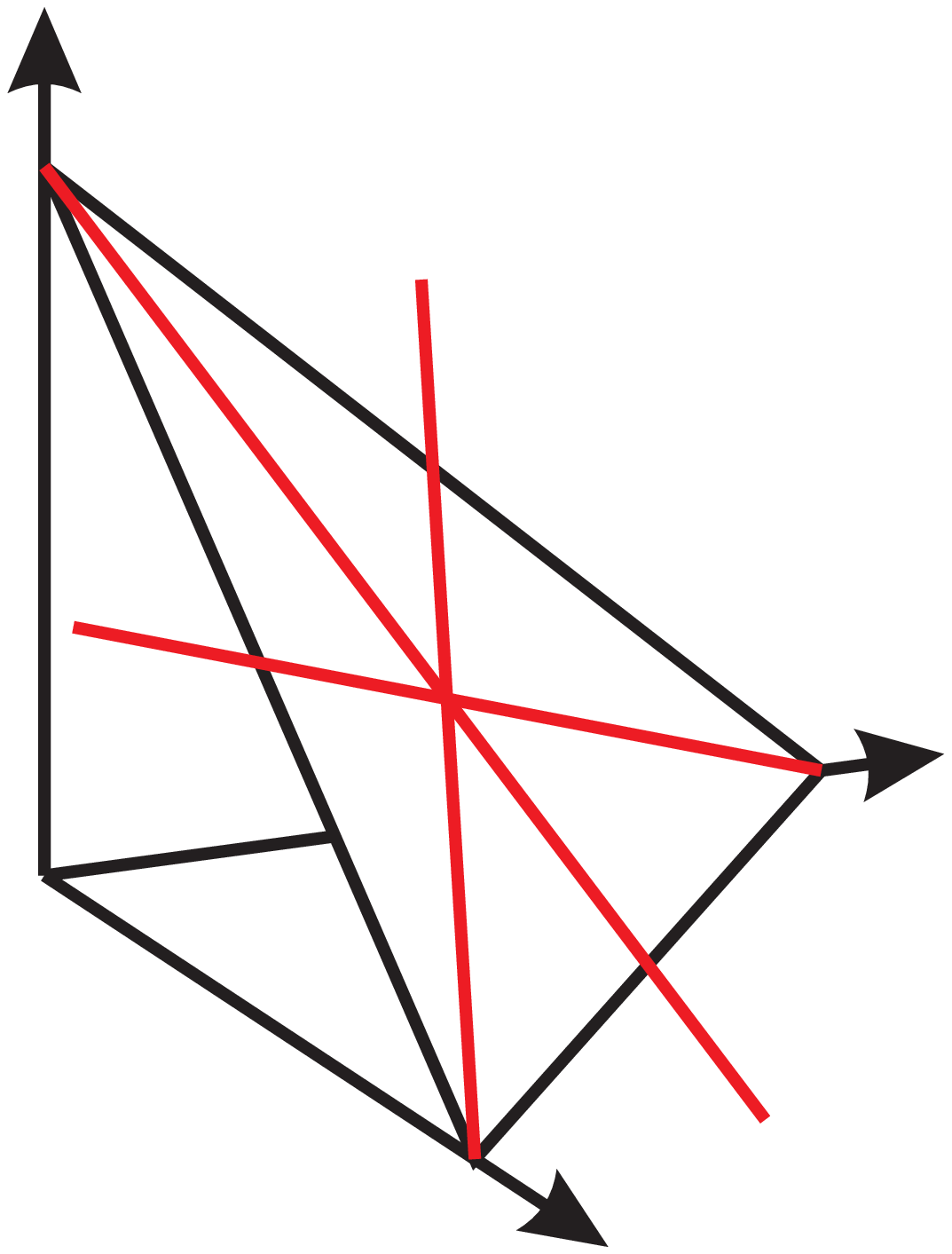 scaled 225}}
}
\rput(10.2,0.4){${\red s_1}$}\rput(9.1,2.9){${\red s_2}$}
}
\end{pspicture}
$$
The element $s_1=(1,2)$ acts as on the left -- as the reflection in the plane with equation
$x_1-x_2=0$. Similarly $s_2=(2,3)$ and $(1,3)$ are reflections in the
planes $x_2-x_3=0$ and $x_1-x_3=0$. These three planes chop the
intersection of
$v^\perp+\frac{1}{3}v$ with the positive quadrant into a
triangle with its boundary barycentrically subdivided (or hexagon). So we 
start to see
the local picture of Figure \ref{fig:localpicture} 
coming from the geometry of these reflecting hyperplanes.

Putting $v^\perp$ into the plane of the page
decomposes the plane into six infinite wedge-shaped regions:
$$
\begin{pspicture}(0,0)(\textwidth,3)
\rput(2.4,1.5){
\rput(-2,-0.8){{\red chambers}}
\rput(0.6,-1.4){${\red s_1}$}\rput(1.95,0.9){${\red s_2}$}
\rput(0,0){\BoxedEPSF{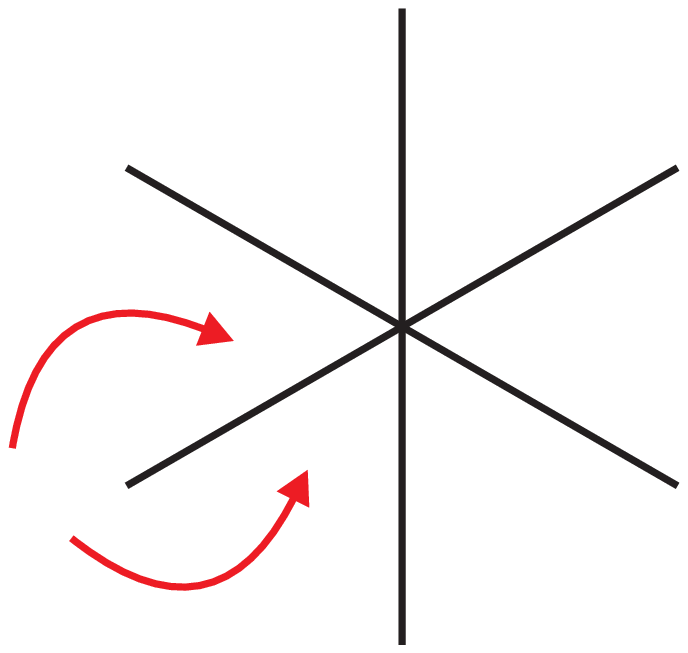 scaled 500}}
}
\rput(-0.5,0){
\psline{->}(4.9,1.5)(6.1,1.5)
\rput(5.5,1.75){intersect}
\rput(5.5,1.25){with $S^1$}
}
\rput(7.25,1.5){
\rput(0,0){\BoxedEPSF{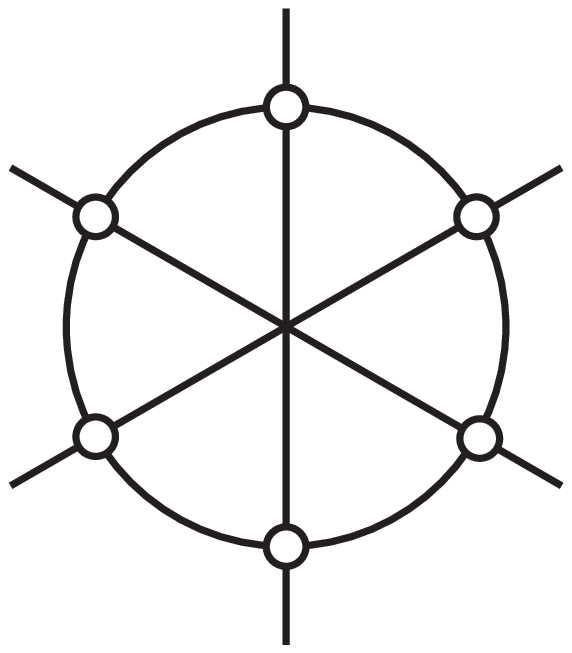 scaled 500}}
}
\rput(9.25,1.5){$\approx$}
\rput(11.25,1.5){
\rput(-1.6,-0.8){{\red chambers}}
\rput(0,0){\BoxedEPSF{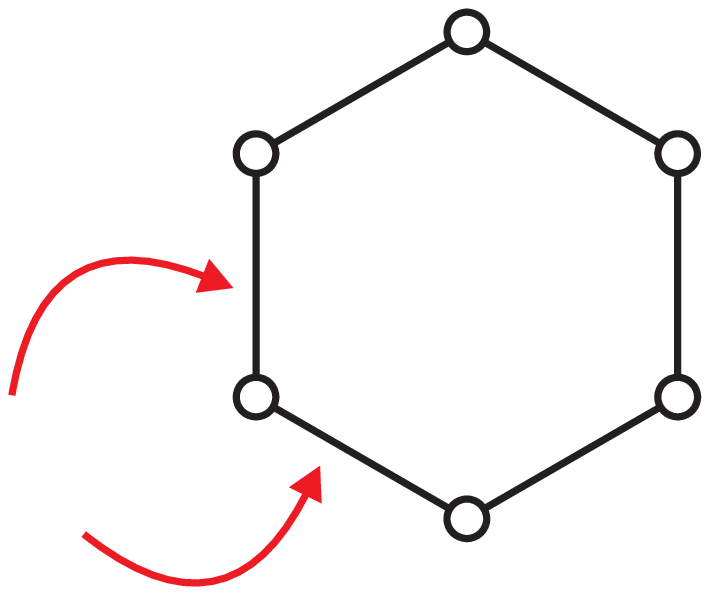 scaled 500}}
}
\end{pspicture}
$$
In the theory of reflection groups (\S \ref{lecture2}) these
regions are also called chambers. The
chambers of our local picture are gotten back by intersecting these
regions with the sphere $S^1$.

\begin{figure}
  \centering
\begin{pspicture}(0,0)(12,5)
\rput(6,4){\BoxedEPSF{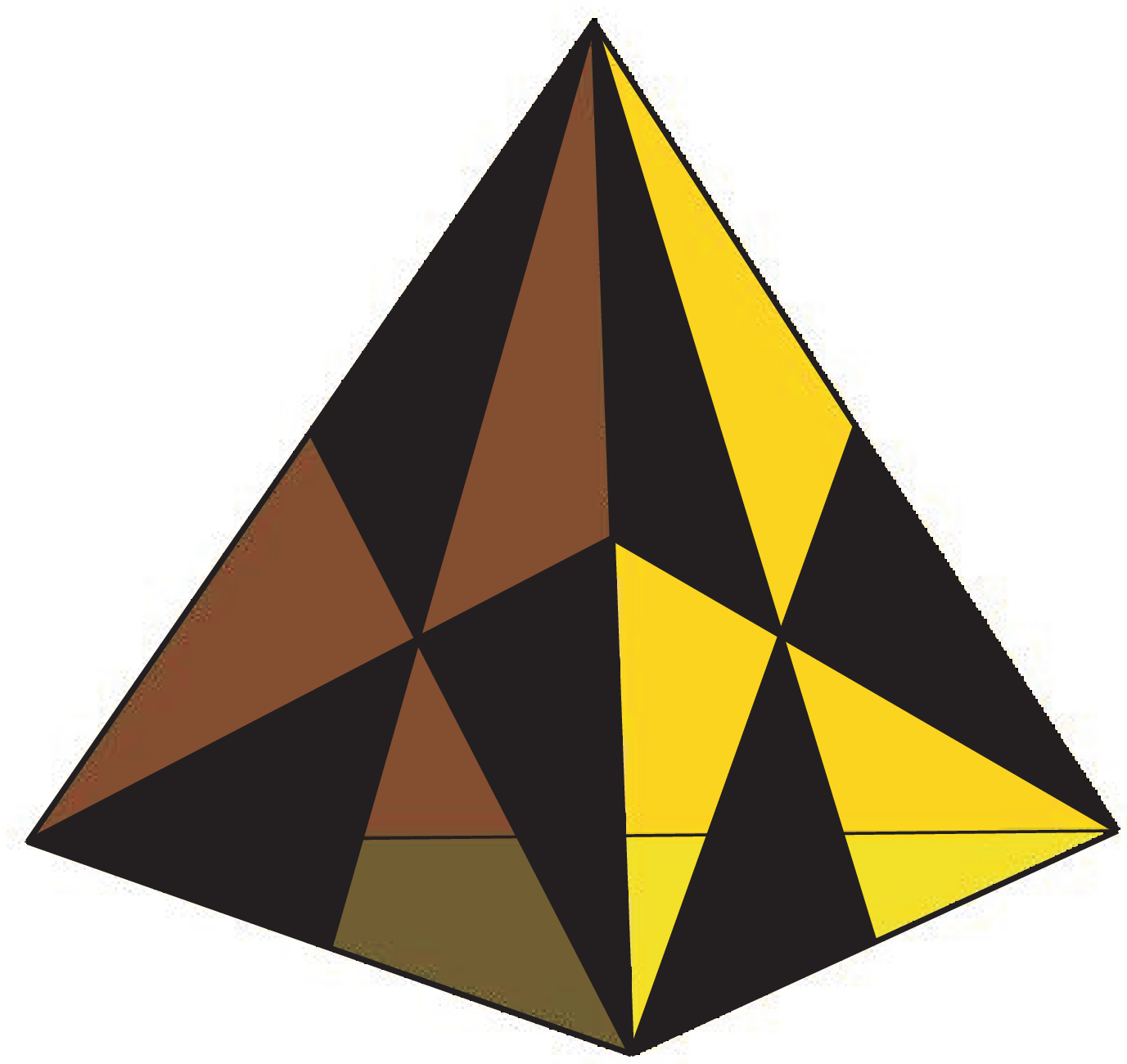 scaled 150}}
\rput(2,2.5){\BoxedEPSF{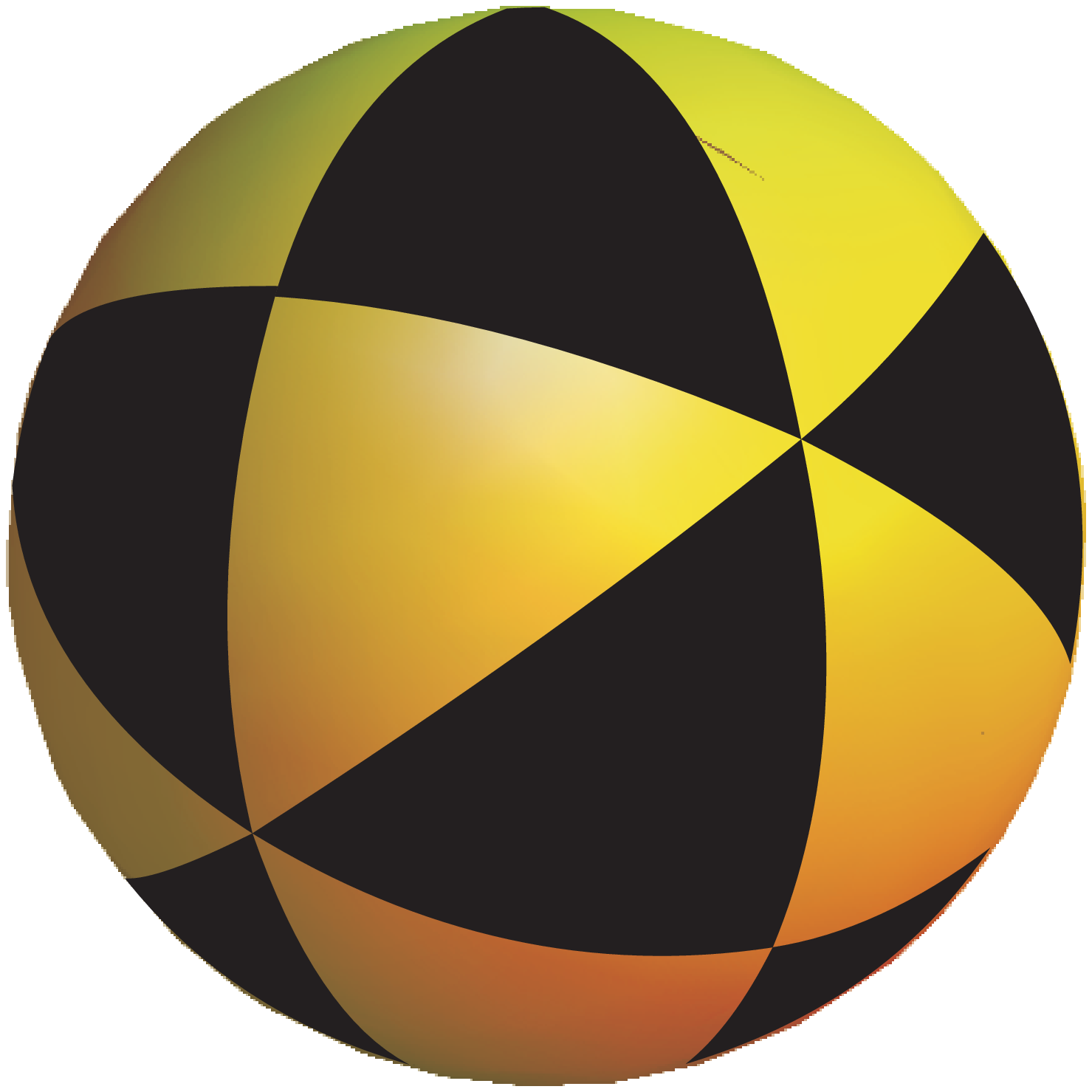 scaled
    300}}
\rput(4.55,4){$\approx$}
\rput(1.3,0){
\rput(8,2.5){\BoxedEPSF{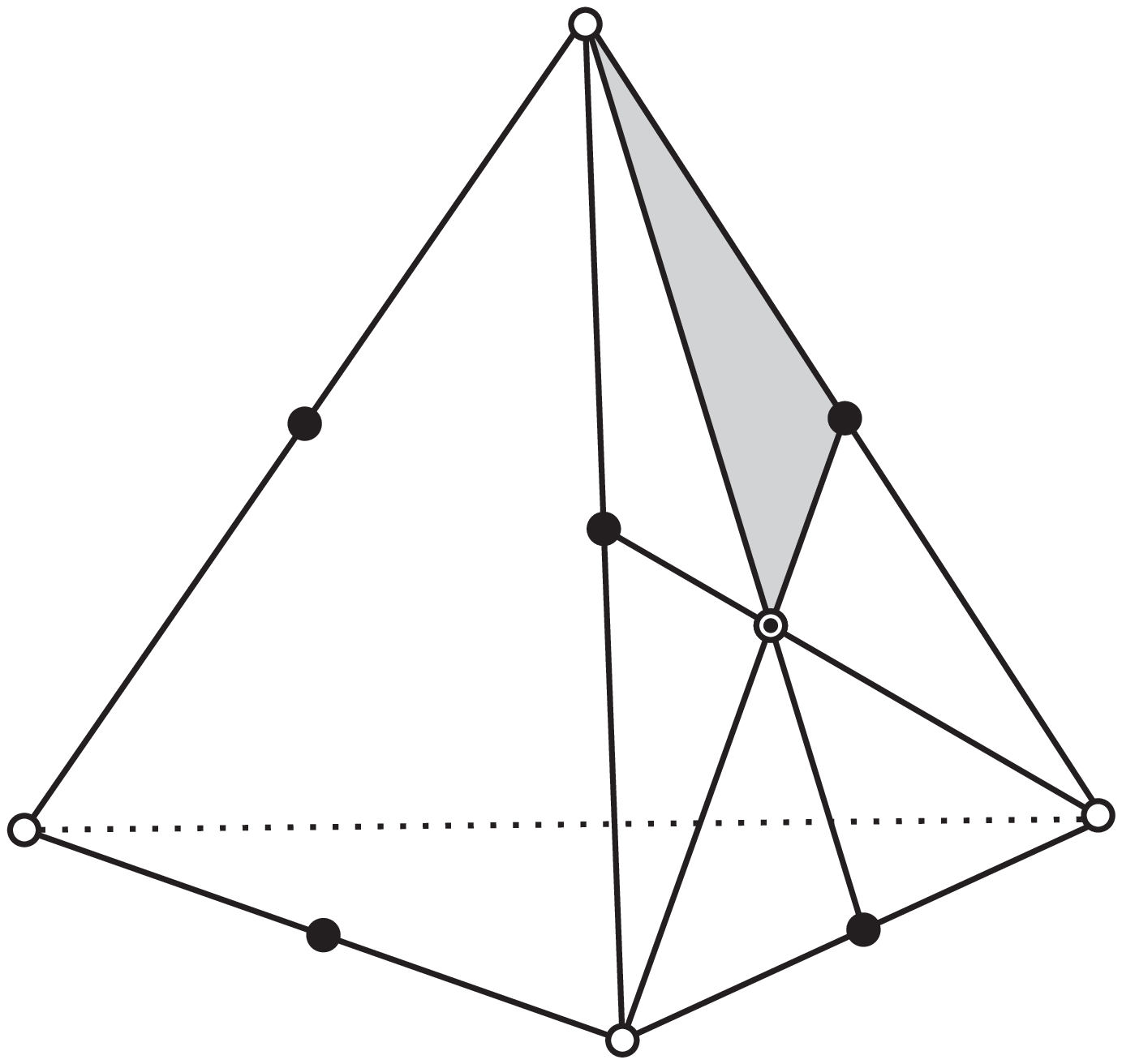 scaled 350}}
\rput(8.1,5){$L_1$}\rput(10.65,1.25){$L_2$}\rput(8.25,0){$L_3$}\rput(5.35,1.25){$L_4$}
\rput(9.95,3){$L_1+L_2$}\rput(9.8,0.5){$L_2+L_3$}\rput(7.5,2.5){$L_1+L_3$}
\rput*(10.1,2.1){$L_1+L_2+L_3$}
}
\end{pspicture}
\caption{The local picture for the flag complex of a four dimensional
  space: the result of intersecting the hyperplanes $x_i-x_j=0$ with
  $S^2$ \emph{(left)}, flattened out \emph{(middle)}, and the picture corresponding
  to the labelled hexagon of Figure \ref{fig:localpicture}
  \emph{(right)}. The shaded $2$-simplex corresponds to the triple
  $L_1\subset L_1+L_2\subset L_1+L_2+L_3$.}
  \label{fig:localpicture2}
\end{figure}

(In the next dimension up we can still draw pictures of some of these
objects. Let $V$ be four dimensional over $k$ and $\Delta$ the two
dimensional simplicial complex with vertices the non-trivial subspaces
of $V$, edges (or $1$-simplicies) the pairs $V_i\subset V_j$ and
$2$-simplicies the triples $V_i\subset V_j\subset V_k$.  We can get
the local picture by working backwards from a symmetric group action
like we did above. If
we have a four dimensional Euclidean space with orthonormal basis
$v_1,v_2,v_3,v_4$, then the convex hull of the $v_i$ is a
tetrahedron lying in the hyperplane $v^\perp+\frac{1}{4}v$ where
$v=v_1+v_2+v_3+v_4$. The six reflecting hyperplanes of the
$\gS_4$-action have equations $x_i-x_j=0$ and slice the boundary of
the tetrahedron barycentrically. Identifying the hyperplane with three
dimensions and intersecting the whole picture with the sphere $S^2$, we
end up with Figure \ref{fig:localpicture2} (left). Flattening it out,
we can retrospectively label the simplicies by lines
$L_1,L_2,L_3,L_4\in V$ and the spaces they generate.)

Returning to the hexagon, the $\gS_3$-action turns out to be regular on the chambers, i.e. given 
chambers $c,c'$ there is a unique $\pi\in\gS_3$ with $\pi c=c'$. This
is most easily seen by brute force: fix a ``fundamental'' chamber $c_0$ and
show that the six elements of $\gS_3$ send it to the six chambers in
the decomposition above. In particular there is a one-one correspondence
between the chambers and the elements of $\gS_3$ given by 
$\pi\in\gS_3\leftrightarrow\mbox{chamber }\pi
c_0$. 

This correspondence gives the adjacency labelings
of the hexagonal local picture of Figure \ref{fig:localpicture}: choose the fixed chamber
$c_0$ to be one
of the two regions bounded by the reflecting lines for
$s_1$ and $s_2$. Starting with the edge of the hexagon contained
in $c_0$, label its vertices by the corresponding reflections as below
left:
$$
\begin{pspicture}(0,0)(\textwidth,3.5)
\rput(1.75,1.75){
\rput(0,0){\BoxedEPSF{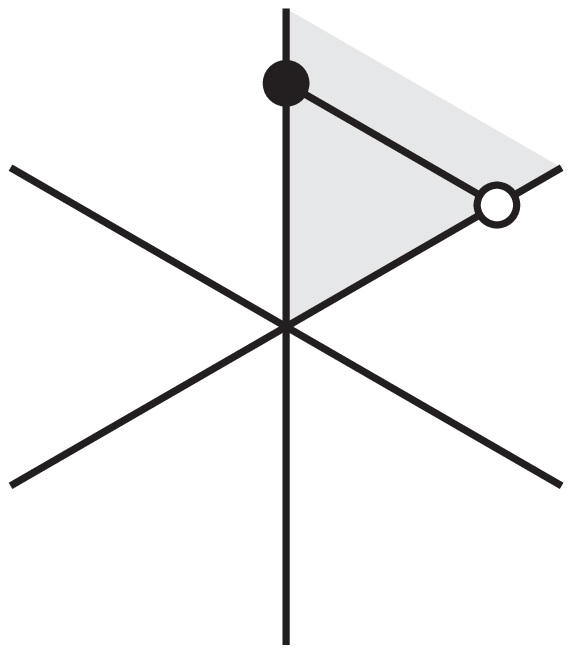 scaled 500}}
\rput(0.75,1.25){$c_0$}
\rput(-0.35,1.25){$s_1$}\rput(1.4,0.5){$s_2$}
}
\rput(-0.25,0){
\rput(6,1.75){
\rput(0,0){\BoxedEPSF{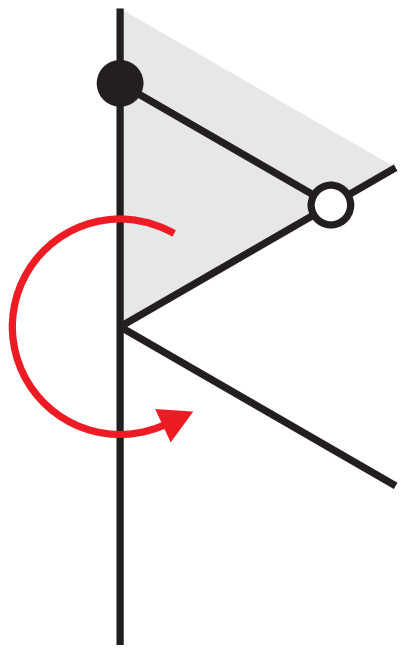 scaled 500}}
\rput(0.4,1.25){$c_0$}
\rput(-0.75,1.25){$s_1$}\rput(1,0.5){$s_2$}
\rput(-1.2,0){${\red \pi}$}
}
\psline{->}(7.5,1.75)(8.8,1.75)
\rput(8.15,2){gives}
\rput(11,1.75){
\rput(0,0){\BoxedEPSF{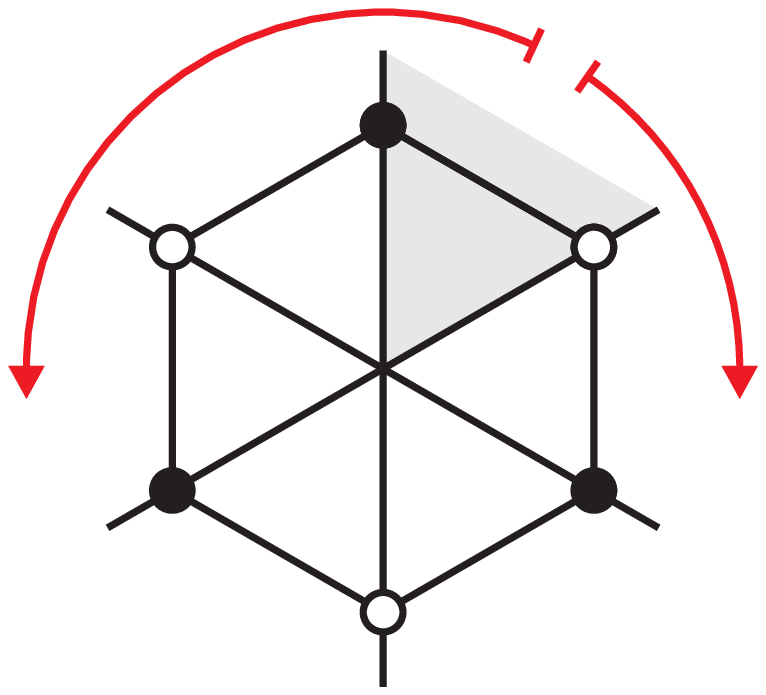 scaled 500}}
\rput(0.3,1.25){$s_1$}\rput(1.45,0.45){$s_2$}
\rput(-1.4,-0.5){$s_1$}\rput(-1.4,0.45){$s_2$}
\rput(-.8,0){$c_1$}\rput(.825,0){$c_2$}
\rput(1.45,-0.5){$s_1$}\rput(0.3,-1.45){$s_2$}
\rput(-1.4,1.5){${\red\pi_1}$}\rput(1.9,0.8){${\red\pi_2}$}
}}
\end{pspicture}
$$
Now transfer this labeled edge to the other chambers using
the $\gS_3$-action as in the picture above middle; the result is shown
above right, where the $i$'s have become $s_i$'s. Vertices on
opposite ends of the same line have different labels
because the antipodal map $x\mapsto -x$ is not
in the action of $\gS_3$ on the plane $v^\perp$.

Finally, to get the metric $\delta$ observe that if $c$ is some
chamber of the local picture in Figure \ref{fig:localpicture}
and $\pi\in\gS_3$ sends $c_0$ to $c$, then $\pi=s_{i_1}\ldots s_{i_k}$
where $s_{i_k},\ldots,s_{i_1}$ are the labels (read from left to right) on the vertices crossed
in a path in the hexagon from $c_0$ to $c$. So for chambers
$c_1,c_2$ we have $\delta(c_1,c_2)=\pi_1^{-1}\pi_2$ where $c_i=\pi_i
c_0$. For our original $c_1,c_2$ we have $\pi_1=s_1s_2$, $\pi_2=s_2$,
hence $\delta(c_1,c_2)=s_2s_1s_2$ as shown in the picture above. 

\paragraph{Second rough definition of building} A building is a set of
chambers with $i$-adjacency, the $i$ coming from some set
$S$, together with a $W$-valued metric $\delta$, for $W$ a reflection
group generated by $S$ and $\delta$ arising from the geometry of $W$.

\paragraph{} 
In the next sections we will make precise and general the
ideas in this rough definition, but working in the reverse order: we
start with reflection groups (\S \ref{lecture2}), then an abstract version of chambers and
adjacency (\S \ref{lecture3}) and finally $W$-valued metrics (\S \ref{lecture4}).

\section{Reflection Groups and Coxeter Groups}
\label{lecture2}

Reflection groups arise as the symmetries of familiar geometric objects;
Coxeter groups are an abstraction of them.
This section covers the basics. All vector spaces and
linear maps here are over the reals $\R$.

A \emph{reflection\/} of a finite dimensional vector space $V$ is a
linear map $s:V\rightarrow V$ for 
which there is a decomposition
\begin{equation}
  \label{eq:1}
V=H_s\oplus L_s  
\end{equation}
where $H_s$ is a hyperplane 
(a codimension one subspace);
$L_s$ is one dimensional; the restriction of $s$ to $H_s$ is the
identity; and the restriction to $L_s$ is the map 
$x\mapsto -x$. Thus a reflection fixes
pointwise a mirror $H_s$, the reflecting hyperplane of $s$,  and acts
as multiplication by $-1$ in some
direction (the reflecting line) not lying in the mirror. 
In particular $s$ is invertible and an involution\footnote{We will
  have no need for them in these notes, but one can reflect a
vector space over an arbitrary field $k$: the definition is identical
except that the restriction of $s$ to the reflecting line $L_s$ is the
map $x\mapsto\zeta x$, where $\zeta$ is a primititve root of unity in
$k$. The only such $\zeta$ in $\R$ is $-1$, hence the definition we have given
of \emph{real\/} reflections. By contrast a complex reflection can
have any finite order.}.

A \emph{reflection group\/} $W$ is a subgroup of $GL(V)$ generated by finitely many
reflections. 

\begin{example}[orthogonal reflections]
\label{example:orthogonal}
The most familiar kind of reflections are the orthogonal
ones for which we further assume that $V$ is a
Euclidean space, i.e. is equipped with an inner product. Then $s$
is orthogonal if in the decomposition (\ref{eq:1}) the line $L_s=H_s^\perp$, 
the orthogonal complement. In particular $L_s$, and hence the 
reflection, is determined by the reflecting
hyperplane, 
unlike a general reflection where both the hyperplane
and the line are needed.
$$
\begin{pspicture}(0,0)(4,3)
\rput(0.5,0){
\rput(1.5,1.5){\BoxedEPSF{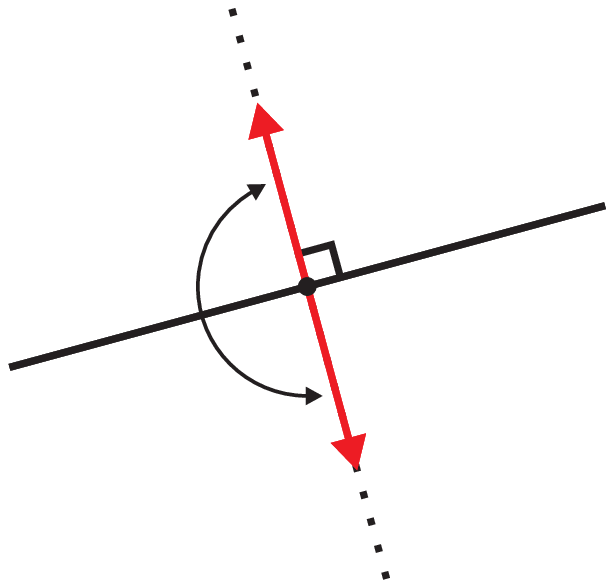 scaled 600}}
\rput(1.5,2.5){$v$}\rput(1.3,0.5){$-v$}
\rput(0.9,2.7){$L_v$}\rput(2.8,1.6){$H_v$}
\rput(0.6,1.6){$s_v$}
}
\end{pspicture}
$$
If $s$ is orthogonal then for any vector $v$ in $L_s$ we have
$s:v\mapsto -v$ with $v^\perp$ fixed pointwise. 
Thus an orthogonal reflection $s$
can be specified by just a non-zero vector $v$, as the reflection with
$H_s=v^\perp$ and $L_s$ spanned by $v$. We write
$s=s_v$, $H_s=H_v$, $L_s=L_v$, and by choosing a sensible basis
one gets that an orthogonal reflection is an orthogonal map of the
Euclidean space.

Let $\HH=\{H_{v_1},\ldots,H_{v_m}\}$ be
hyperplanes in Euclidean $V$ and $W$ the reflection group generated by the orthogonal
reflections $s_{v_1},\ldots,s_{v_m}$. As an exercise the reader can show that if $W\HH=\HH$, i.e. $g
H_{v_i}=H_{v_j}$ for all $g\in W$ and all $v_i$,
then $W$ is finite (\emph{hint\/}: $|W|\leq (2m)!$). It turns out 
(although this is harder) that 
$\HH$ then consists of \emph{all\/} the reflecting hyperplanes of $W$.
\end{example}

\begin{example}[a finite reflection group]
\label{example:finite}
Let
$V$ be Euclidean with an orthonormal basis $v_1,\ldots,v_{n+1}$ and
$\HH$ the hyperplanes $H_{ij}:=(v_i-v_j)^\perp$
for $1\leq i\not=j\leq n+1$ (in other words, $H_{ij}$ is the hyperplane with
equation $x_i-x_j=0$). The reflection $s_{v_i-v_j}$ sends
$v_i-v_j$ to $v_j-v_i$, thus swapping the vectors $v_i$
and $v_j$. Any other basis vector is orthogonal to $v_i-v_j$, so lies
in $H_{ij}$, and  is fixed. Thus if $\pi=(i,j)\in\gS_{n+1}$ 
then $s_{v_i-v_j}H_{k\ell}=H_{\pi(k),\pi(\ell)}$. 

Now let $W$ be the group generated by the reflections
$s_{v_i-v_j}$. We have just shown that $W\HH=\HH$, so $W$ is a
finite reflection group by the exercise above. Indeed, $W$ is the symmetric group
$\gS_{n+1}$ acting by permuting coordinates as in 
\S \ref{lecture1}. To make this identification we have already seen that each
$s_{v_i-v_j}$, and so every element of $W$, permutes 
the basis vectors $v_1,\ldots,v_{n+1}$. This gives a
homomorphism $W\rightarrow\gS_{n+1}$. Injectivity of this homomorphism
follows as the
$v_i$ span $V$ and surjectivity as the transpositions $(i,j)$ generate
$\gS_{n+1}$. 

The convex hull of the $v_i$ is the standard $n$-simplex,
barycentrically subdivided by its $n(n-1)$ hyperplanes of reflectional
symmetry (the $\HH$),
each of which
is a reflecting hyperplane of $W$. This is the picture we had for
$n=2$ and $n=3$ in \S \ref{lecture1}. Finite reflection groups are often
called \emph{spherical\/} as the geometrical realisation of their
Coxeter complexes (the boundary of 
the barycentrically divided $n$-simplex in this case; see Example 
\ref{example:coxeter.complexes}
for the general definition) are spheres.
\end{example}

\begin{example}[an affine reflection group] 
\label{example:affine}
Let $V$ be $2$-dimensional and consider reflections $s_0,s_1$ where the reflecting
hyperplanes and lines are shown below left (there is no inner product
this time).
The reflecting hyperplanes are different but both have the same reflecting line:
$L_{s_0}=L=L_{s_1}$. If $W$ is the group generated by $s_0,s_1$ then
$W$ leaves invariant any affine line parallel to $L$
as the $s_i$ do. But if $\HH=\{H_{s_0},H_{s_1}\}$ then
$W\HH\not=\HH$ as $s_0H_{s_1}\not\in\HH$. Indeed, we must expand $\HH$
to the infinite set shown below right before it becomes
$W$-invariant:
$$
\begin{pspicture}(0,0)(\textwidth,3.5)
\rput(-1,0){
\rput(3,1.5){\BoxedEPSF{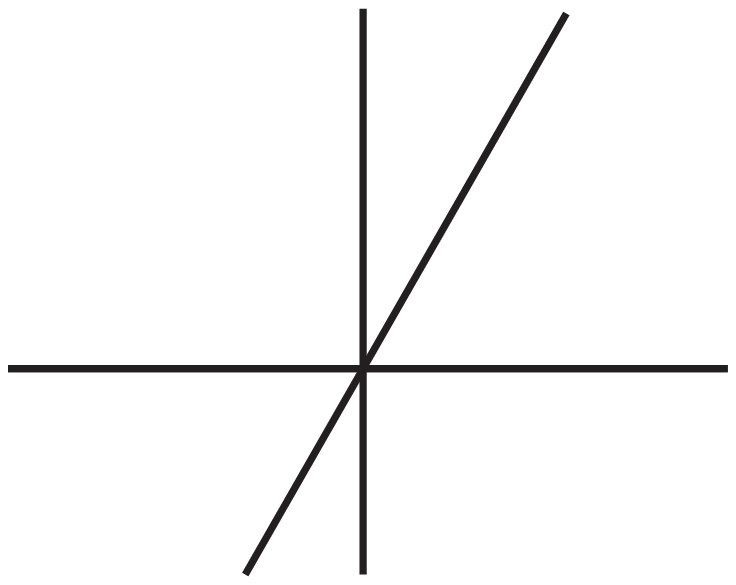 scaled 500}}
\rput(2.6,2.8){$H_{s_0}$}\rput(4.3,2.8){$H_{s_1}$}
\rput(4.5,1.35){$L_{s_0}=L_{s_1}$}
\rput(4.5,0.25){$W=\langle s_0,s_1\rangle$}
}
\rput(5.5,0){
\rput(3,1.5){\BoxedEPSF{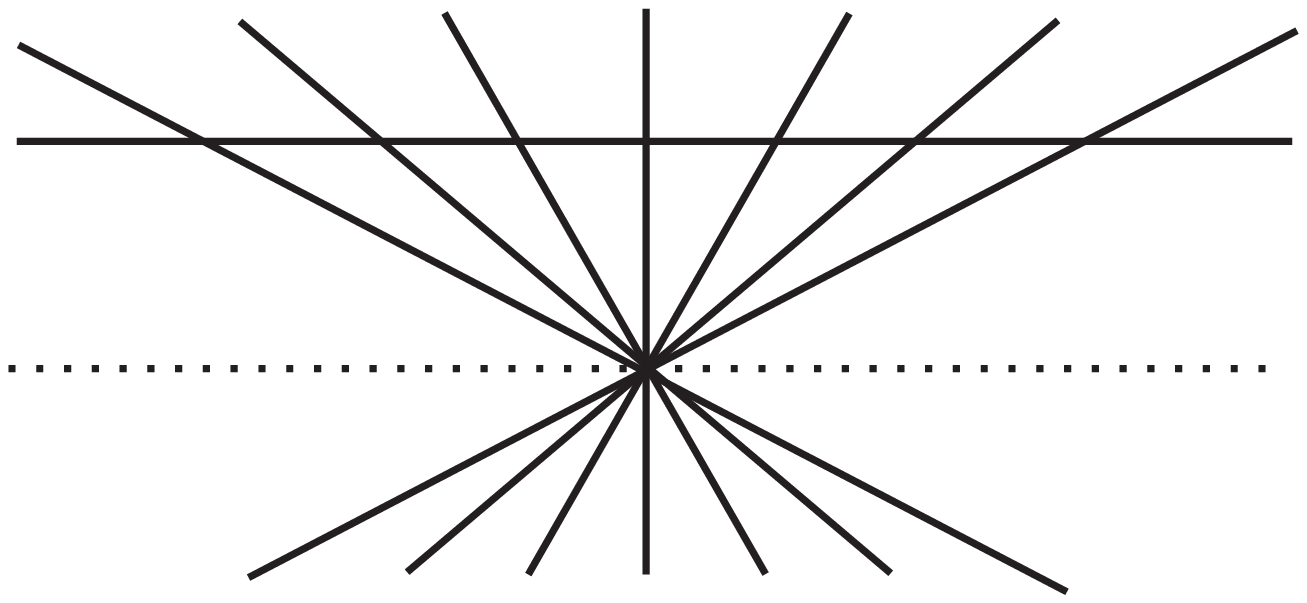 scaled 500}}
}
\rput(-1.5,0){
\rput(10,3.2){$s_0$}\rput(11,3.2){$s_1$}
\rput(9,3.2){$s_0s_1s_0$}\rput(12,3.2){$s_1s_0s_1$}
\rput(8,3.2){$\ldots$}\rput(13,3.2){$\ldots$}
\rput(13,1.5){invariant affine line}
\psline{->}(13,1.65)(13,2.3)
\rput(7.5,1.4){$L$}
}
\end{pspicture}
$$
In fact, by identifying the invariant affine line with 
the reals, $W$ is isomorphic to the
group of ``affine reflections'' of $\R$ in the integers $\Z$, 
i.e. to
the group of transformations of $\R$ generated by the maps $s_n:x\mapsto 2n-x$
for $n\in\Z$. 
The element $s_1s_0$ acts on the affine line as the
translation $x\mapsto x+2$ so has infinite order. In particular $W$ is
infinite. This also follows from $\HH$ being infinite as the
reflections in the hyperplanes in $\HH$ are the $W$-conjugates of $s_0,s_1$.
\end{example}

\begin{example}[hyperbolic reflections]
\label{hyperbolic.reflection}
Let $V$ be $3$-dimensional and again there is no inner 
product. 
Let
$a,b,c$ be real numbers such that $a^2+b^2>c^2$, and consider the
reflection $s$ with reflecting hyperplane $H_s$ having the equation
$ax+by-cz=0$ and reflecting line $L_s$ spanned by the vector 
$v=(a,b,c)$. 
Then $v$ lies on the outside of the pair of cones with equation 
$z^2=x^2+y^2$ and $H_s$ passes through the interior of this cone:
$$
\begin{pspicture}(0,0)(\textwidth,4)
\rput(2.5,1.8){
\rput(0,0){\BoxedEPSF{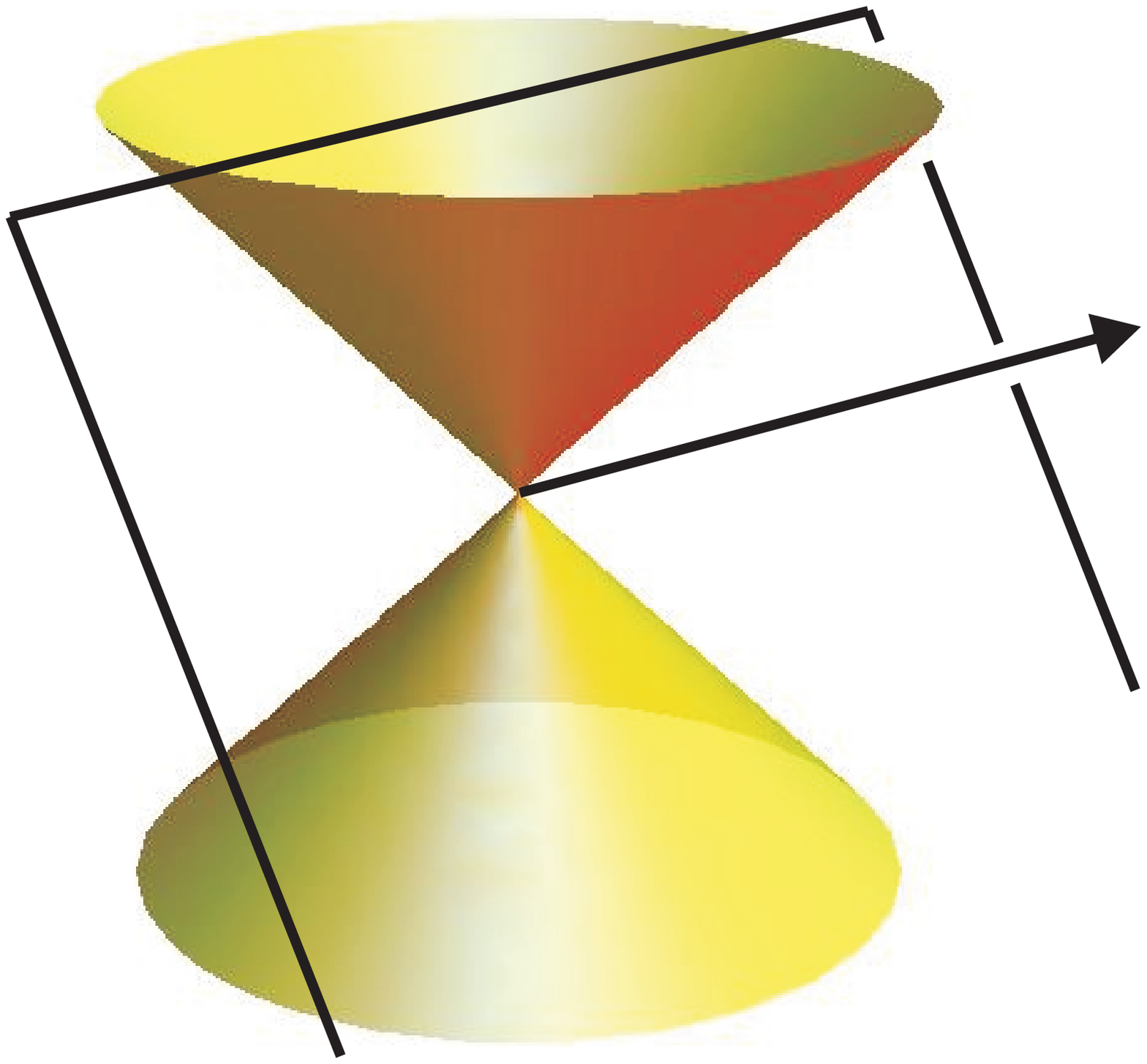 scaled 250}}
\rput(1.5,-1.8){$z^2=x^2+y^2$}
\rput(-2,0.5){$H_s$}\rput(2.5,1){$v$}
}
\rput(8.5,1.8){
\rput(0,0){\BoxedEPSF{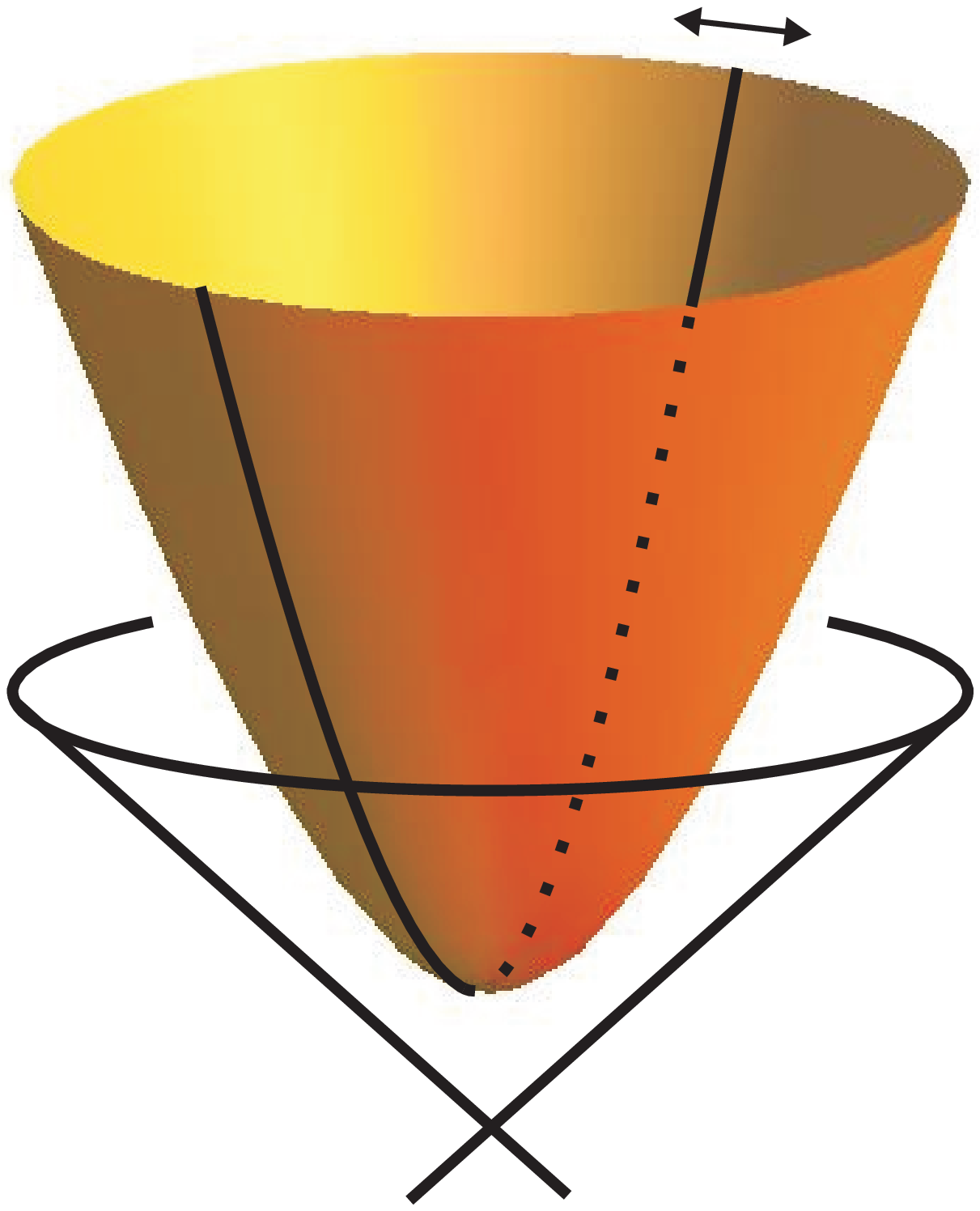 scaled 250}}
\rput(0.8,2.3){$s$}
\rput(2.8,1){$x^2+y^2-z^2=-1$}
}
\end{pspicture}
$$
One can check that $s$ leaves invariant each sheet of the  
two sheeted hyperboloid with equation
$x^2+y^2-z^2=-1$. 
Either sheet is a model for the hyperbolic plane.
Intersecting $H_s$ with the top sheet gives a
hyperbola -- a straight line of hyperbolic geometry -- and $s$
is the ``hyperbolic reflection" of the plane in this 
line\footnote{Although there is no inner product in Examples 
\ref{example:affine} and \ref{hyperbolic.reflection}, it is possible
to endow $V$ with a bilinear form so that the reflections are 
``orthogonal" with respect to this form.}.
\end{example}

Returning to the finite orthogonal case, let $V$ be Euclidean,
$\HH=\{H_i\}_{i\in T}$ a finite set of hyperplanes and $W=\langle
s_i\rangle_{i\in T}$ the group
generated by the orthogonal reflections in the $H_i$. Suppose also that
$W\HH=\HH$, so $W$ is finite and $\HH$ is the set of all
reflecting hyperplanes of $W$ as above. 

For each $i\in T$ choose a linear functional $\aa_i\in V^*$ with
$H_i=\ker\aa_i$. The choice of $\aa_i$ is unique upto scalar
multiple and $H_i$ consists of those $v\in V$ with $\aa_i(v)=0$.
The two sides (or half-spaces) of the hyperplane consist of the
$v$ with $\aa_i(v)>0$ or the $v$ with $\aa_i(v)<0$. 

Fix an $T$-tuple $\ve=(\ve_i)_{i\in T}$, with $\ve_i\in\{\pm
1\}$, and consider the set
\begin{equation}
  \label{eq:2}
c=c(\ve)=\{v\in V\,|\,
\ve_i\aa_i(v)>0\text{ for all }i
\}.
\end{equation}
So each $\aa_i(v)$ is non-zero and $\aa_i(v),\ve_i$ have the same sign
for all $i$.
If this set is non-empty then call it
a \emph{chamber\/} of $W$. A non-empty set of the form
\begin{equation}
  \label{eq:3}
a=a(\ve)=\{v\in V\,|\,\aa_{i_0}(v)=0\text{ for some }i_0,
\text{ and $\ve_i\aa_i(v)>0$ for all $i\not=i_0$}
\}  
\end{equation}
is called a \emph{panel\/}. Here is the example from \S \ref{lecture1}:
$$
\begin{pspicture}(0,0)(\textwidth,4)
\rput(3,2){
\rput(0,0){\BoxedEPSF{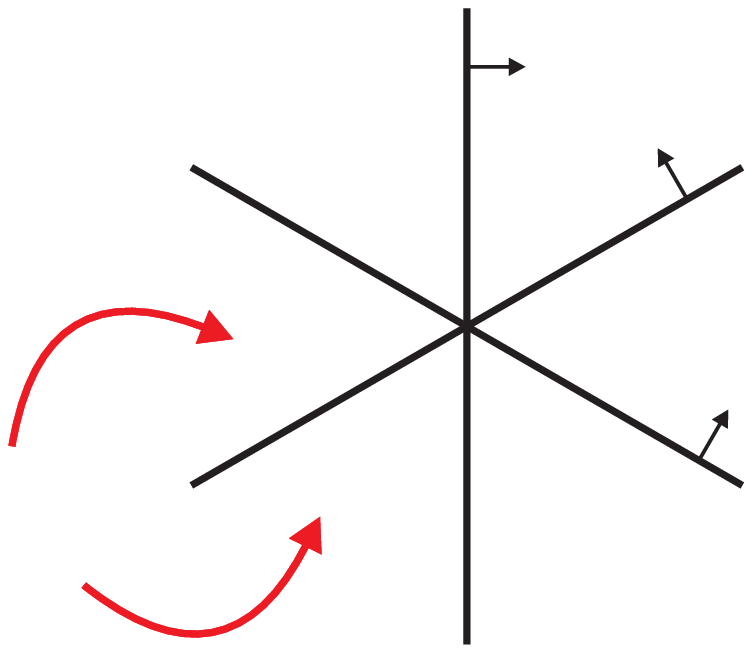 scaled 600}}
}
\rput(-3.3,0){
\rput(4.25,0.8){{\red chambers}}
\rput(7.4,3){$\scriptstyle{+++}$}
\rput(8,2){$\scriptstyle{+-+}$}
\rput(7.4,1){$\scriptstyle{+--}$}
\rput(6.2,3){$\scriptstyle{-++}$}
\rput(5.8,2){$\scriptstyle{-+-}$}
\rput(6.2,1){$\scriptstyle{---}$}
\rput(7.5,3.5){$\aa_1$}\rput(8.35,3.2){$\aa_2$}\rput(8.7,1.3){$\aa_3$}
}
\rput(-1,0){
\rput(10.5,2){
\rput(0,0){\BoxedEPSF{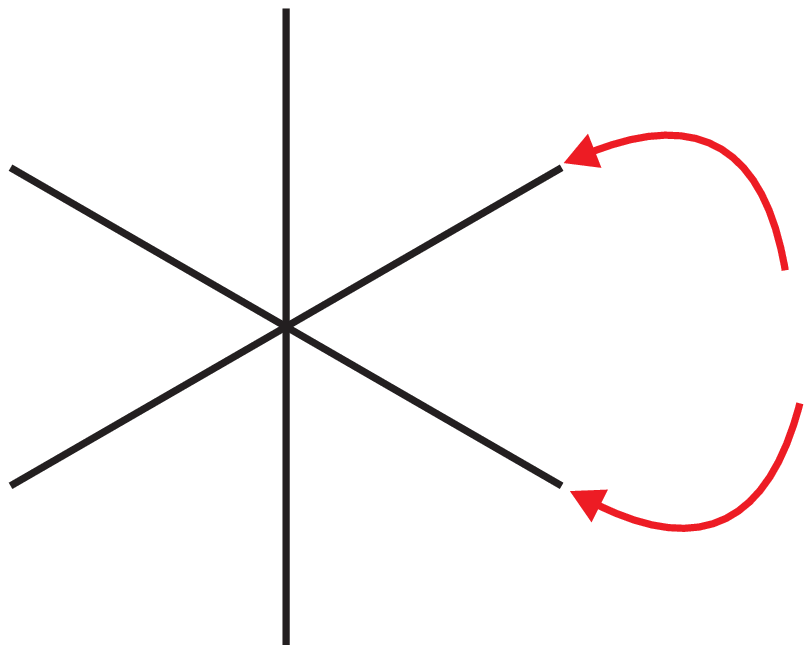 scaled 600}}
}
\rput(3,0){
\rput*(6.8,0.75){$0--$}
\rput*(6.8,3.25){$0++$}
\rput*(5.6,1.4){$-0-$}
\rput*(8,1.4){$+-0$}
\rput*(5.6,2.6){$-+0$}
\rput*(8,2.6){$+0+$}
\rput(10,2){{\red panels}}
}}
\end{pspicture}
$$
where there are three hyperplanes in $\HH$ and the $\aa_i$ are chosen
so that $\aa_i(v)>0$ for those $v$ on the side indicated by the
arrow. The chambers are marked by their $T$-tuples. 
There are $2^3$ $T$-tuples but only $6$ chambers
because the tuples $++-$ and $--+$ give empty sets in (\ref{eq:2}).
Extend the notation to include panels (\ref{eq:3}) by placing a $0$ in
the $i_0$-th position. There are then $3.2^2$ such tuples but only
$6$ give non-empty panels, with two lying on each reflecting line.

There is an obvious notion of adjacency between chambers 
suggested by these pictures. Say that $a$ is a panel
of the chamber $c$ if the corresponding $T$-tuples are identical
except in one position where the tuple for $a$ has a $0$. It
turns out that this can also be defined topologically: $a$ is a panel
of $c$ exactly when $\bar{a}\subset\bar{c}$, the closures of these
sets with respect to the usual topology on $V$. 

Chambers $c_1$ and $c_2$ are then \emph{adjacent\/} if they
share a common panel. In the Example from \S \ref{lecture1}, chambers are adjacent
when they share a common edge. 

The
adjacency relation can be refined by bringing the 
reflection group $W$ into the picture. In \S \ref{lecture1} we saw
that $\gS_3$ acts regularly as a reflection group
on the chambers. This turns out to be true in general
for the $W$-action on the chambers:
given chambers $c,c'$ there is a unique $g\in W$ with $g c=c'$. 
Fix one of the chambers $c_0$. Then the regular action gives
the chambers are in one-one correspondence with the
elements of $W$ via $g\in W\leftrightarrow \text{chamber }gc_0$.

Now let
$S=\{s_1,\ldots,s_n\}$ be those reflections in $W$ whose
hyperplanes $H_1,\ldots,H_n$ are spanned by a panel of the fixed chamber $c_0$. Thus
$S=\{s_1,s_2\}$ for the $c_0$ in the example from \S \ref{lecture1}:
$$
\begin{pspicture}(0,0)(\textwidth,4)
\rput(3,2){
\rput(0,0){\BoxedEPSF{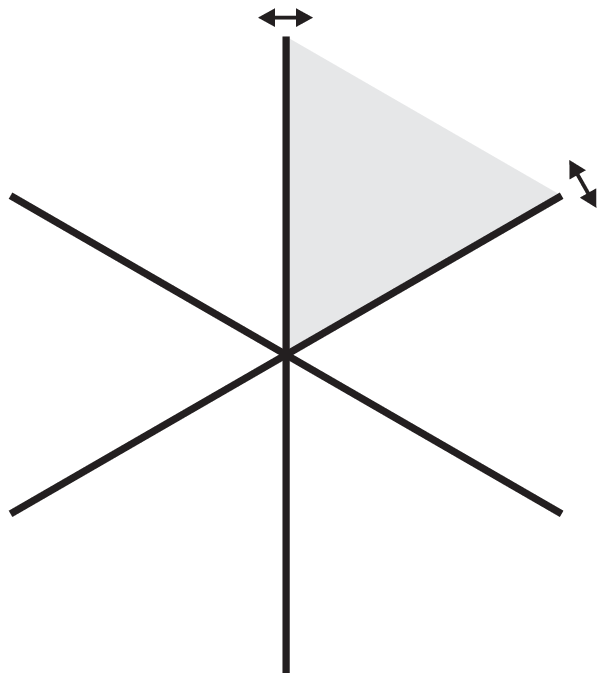 scaled 600}}
\rput(0.5,1){$c_0$}
\rput(-0.5,2){$s_1$}\rput(2,1){$s_2$}
}
\rput(9.5,2){
\rput(0,0){\BoxedEPSF{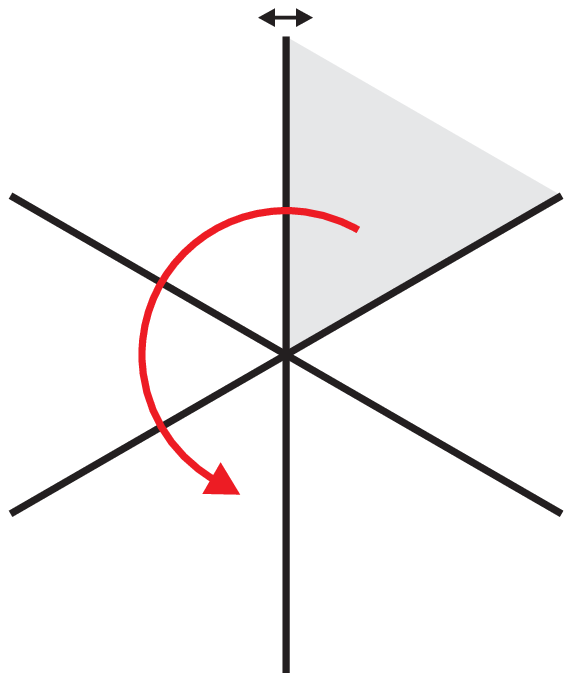 scaled 600}}
\rput(0.7,1){$c_0$}\rput(-0.7,1){$sc_0$}
\rput(-0.7,-1){$c_1$}\rput(0.7,-1){$c_2$}
\rput(-0.4,2){$s$}
\rput(-1.1,0){${\red g}$}
}
\end{pspicture}
$$
Suppose $c_1,c_2$ are a pair of adjacent chambers as 
above right. Then there is a $g\in W$ with 
$c_1=gc_0$.
Translating the picture back to $c_0$ we have $g^{-1}c_1=c_0$
and $g^{-1}c_2$ are adjacent chambers, and the common panel of $c_1,c_2$
is sent by $g^{-1}$ to a common panel of $c_0$ and $g^{-1}c_2$
(these are most easily seen using the topological version of
adjacency). If $s\in S$ is the reflection in the hyperplane spanned by
the common panel of $c_0$ and $g^{-1}c_2$, then the chamber $g^{-1}c_2$ is the same as the
chamber $sc_0$. 

Thus $c_1=gc_0$, $c_2=(gs)c_0$, and we have the following more refined
description of adjancey:
\begin{equation}
  \label{eq:4}
\text{the chambers adjacent to the chamber $gc_0$ are the
  $(gs)c_0$ for
  $s\in S$.}
\end{equation}
When $S=\{s_1,\ldots,s_n\}$ we say that chambers $gc_0$ and $gs_ic_0$ are $i$-adjacent.
In our running example, the chambers adjacent to $gc_0$ are $gs_1c_0$
and $gs_2c_0$, and these are the two that were $1$- and $2$-adjacent
to $gc_0$ in \S \ref{lecture1}. 

\paragraph{Coxeter groups} We motivate the definition of Coxeter
group by quoting two
facts, staying with the assumptions above where $W$ is generated
by orthogonal reflections in finitely many hyperplanes $\HH$ with $W\HH=\HH$:

\paragraph{Fact 1.}  The group $W$ is generated by the reflections
$s\in S$ in those hyperplanes spanned by a panel of the fixed chamber $c_0$.

\paragraph{}
In our running example we can see a
how a proof might work using induction on the ``distance'' of
a chamber from $c_0$.  If $g$ is an
element of $W$ then there is a
chamber adjacent to the chamber $gc_0$ that is closer to 
$c_0$ than $gc_0$ is. If this closer chamber is $g'c_0$ say, then by
(\ref{eq:4}) we have $g=g's$ for some $s\in S$. Repeat the
process until $g$ completely decomposes
as a word in the $s\in S$. 

\paragraph{Fact 2.} With respect to the generators $S$ 
the group $W$ admits a presentation 
\begin{equation}
  \label{eq:5}
\langle s\in S\,|\, (s_i s_j)^{m_{ij}}=1\rangle  
\end{equation}
where the $m_{ij}\in\Z^{\geq 1}$ and are such that
$m_{ij}=m_{ji}$, and $m_{ij}=1$ if and only if $i=j$
(so in particular, $s_i^2=1$).

\paragraph{}
If $s_i$ and $s_j$ are reflections in $W$ finite, then
the element $s_i s_j$ has finite order $m_{ij}\geq
2$. 
So the relations in the presentation (\ref{eq:5}) certainly hold. The
content of Fact 2 is that these relations suffice.
Geometrically, $s_i s_j$ is a rotation ``about'' the
intersection $H_i\cap H_j$ of the corresponding hyperplanes. 


In Example \ref{example:finite} we have $S=\{s_1,\ldots,s_n\}$ where
$s_i=s_{v_i-v_{i+1}}$. The $s_is_{i+1}$ have order $3$ and all other
$s_is_j$ have order $2$. Moreover $W$ is isomorphic to $\gS_{n+1}$ via
the map induced by $(i,i+1)\mapsto s_i$. Our running example of the
action of $\gS_3$ on $3$-dimensional $V$ is the $n=2$ case of this. 

Here is the promised abstraction of reflection group: a group $W$ is
called a \emph{Coxeter group\/} if it admits a presentation
(\ref{eq:5}) with respect to some finite $S$, where the
$m_{ij}\in\Z^{\geq 1}\cup\{\infty\}$ satisfy the rules following
(\ref{eq:5}). 
Sometimes the dependency on the relations $S$ is emphasized and 
$(W,S)$ is called a Coxeter \emph{system\/}.

We want the new concept to cover all the examples we have seen so far
in this section, including the affine group in Example
\ref{example:affine} where the element $s_1s_0$ had infinite
order. This is why in the definition of Coxeter group
the conditions on the $m_{ij}$ are relaxed to allow them to be
infinite. A 
relation $(s_i s_j)^{m_{ij}}=1$
is omitted from the presentation when $m_{ij}=\infty$. 

There is a standard shorthand for a Coxeter presentation
(\ref{eq:5}) called the Coxeter symbol. This is a graph
with nodes the $s_i\in S$, and where nodes $s_i$ and $s_j$ are joined by an edge
labeled $m_{ij}$ if $m_{ij}\geq 4$, an unlabeled edge if 
$m_{ij}=3$ and no edge when $m_{ij}=2$:
$$
\begin{pspicture}(0,0)(1,1)
\pscircle[fillstyle=solid,fillcolor=white](0,0.8){0.125}
\pscircle[fillstyle=solid,fillcolor=white](1,0.8){0.125}
\rput(0.5,0.2){$m_{ij}=2$}
\end{pspicture}
\quad\quad\quad
\begin{pspicture}(0,0)(1,1)
\psline(0,0.8)(1,0.8)
\pscircle[fillstyle=solid,fillcolor=white](0,0.8){0.125}
\pscircle[fillstyle=solid,fillcolor=white](1,0.8){0.125}
\rput(0.5,0.2){$m_{ij}=3$}
\end{pspicture}
\quad\quad\quad
\begin{pspicture}(0,0)(1,1)
\psline(0,0.8)(1,0.8)
\pscircle[fillstyle=solid,fillcolor=white](0,0.8){0.125}
\pscircle[fillstyle=solid,fillcolor=white](1,0.8){0.125}
\rput(0.5,1){$m_{ij}$}
\rput(0.5,0.2){$m_{ij}\geq 4$}
\end{pspicture}
$$
The examples from \S \ref{lecture1} and Example \ref{example:affine}
are then:
$$
\begin{pspicture}(0,0)(\textwidth,4)
\rput(3.5,2){
\rput(0,0){\BoxedEPSF{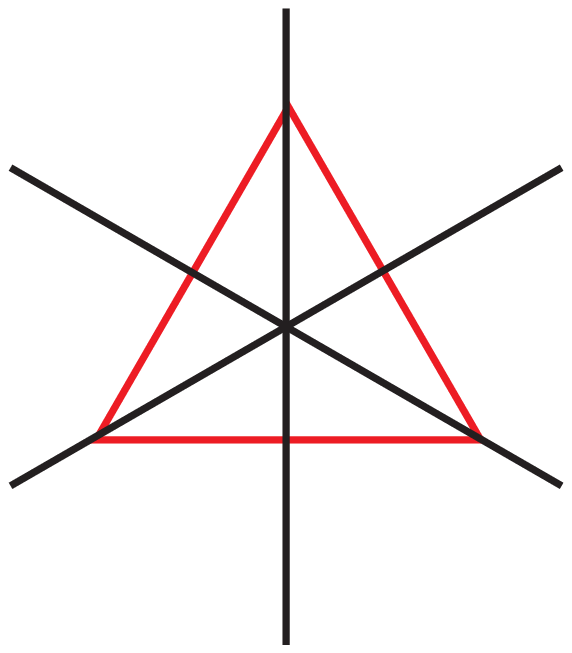 scaled 600}}
\rput(-0.3,1.85){$s_1$}\rput(2,1){$s_2$}
}
\rput(0.5,2){
\psline(0,0)(1,0)
\pscircle[fillstyle=solid,fillcolor=white](0,0){0.125}
\pscircle[fillstyle=solid,fillcolor=white](1,0){0.125}
\rput(0,-0.3){$s_1$}\rput(1,-0.3){$s_2$}
}
\rput(11.5,1){
\psline(0,0)(1,0)
\pscircle[fillstyle=solid,fillcolor=white](0,0){0.125}
\pscircle[fillstyle=solid,fillcolor=white](1,0){0.125}
\rput(0,-0.3){$s_0$}\rput(1,-0.3){$s_1$}
\rput(0.5,0.2){$\infty$}
}
\rput(9.5,2){
\rput(0,0){\BoxedEPSF{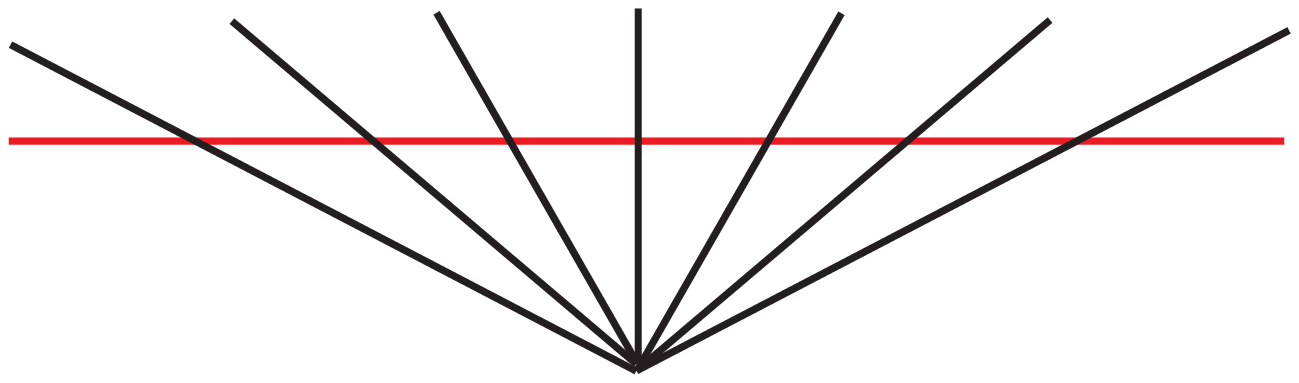 scaled 500}}
\rput(0,1.2){$s_0$}\rput(1,1.2){$s_1$}
}
\end{pspicture}
$$

\begin{remark}\label{remark:titsrep}
What is the relationship between the concrete reflection groups defined at the
beginning of this section and the abstract Coxeter groups defined at
the end? The answer is that the
Coxeter groups are \emph{discrete\/} reflection groups: for a Coxeter system $(W,S)$
one can construct a faithful representation $(W,S)\rightarrow GL(V)$
for some vector space $V$, where the $s\in S$ act on $V$ as
reflections, and the image of $(W,S)$ is a discrete subgroup of
$GL(V)$. 
\end{remark}

\section{Chamber Systems and Coxeter Complexes}
\label{lecture3}

We have seen several examples of sets of chambers with
different kinds of adjacency between them. 
This section introduces the formalization of this idea: chamber systems.

A \emph{chamber system\/} over a finite set $I$ is a set $\Delta$ equipped with
equivalence relations $\sim_i$, one for each $i\in I$. The
$c\in\Delta$ are the \emph{chambers\/} and two chambers are \emph{$i$-adjacent\/} when
$c\sim_i c'$.

The generic picture to keep in mind is below where
chambers are $i$-adjacent if they share a common $i$-labeled edge.
Thus, $c_0\sim_1 c_1,c_0\sim_2 c_2$, etc.
$$
\begin{pspicture}(0,0)(\textwidth,3)
\rput(6.25,1.5){
\rput(0,0){\BoxedEPSF{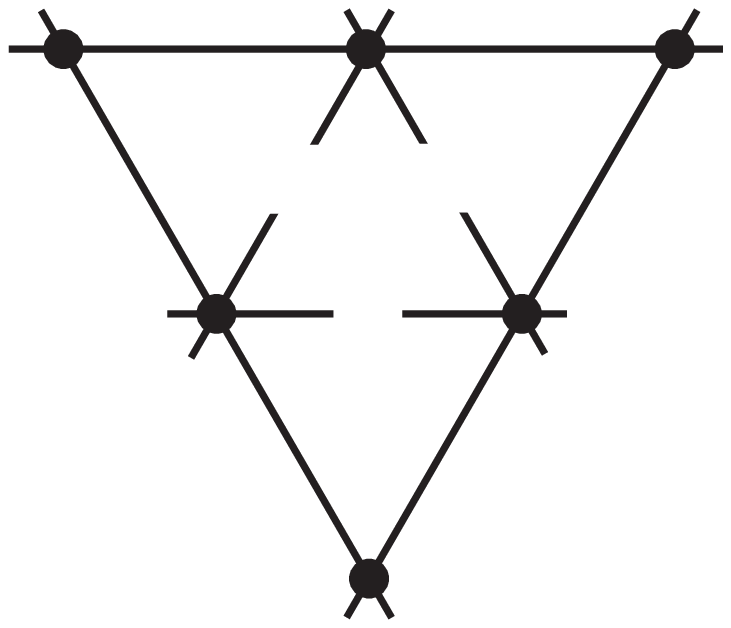 scaled 500}}
\rput(0.425,0.675){${\scriptstyle 1}$}
\rput(-0.375,0.675){${\scriptstyle 2}$}
\rput(0,0){${\scriptstyle 3}$}
\rput(0,0.4){${\scriptstyle c_0}$}
\rput(0.75,0.8){${\scriptstyle c_1}$}
\rput(-0.75,0.8){${\scriptstyle c_2}$}
\rput(0,-0.5){${\scriptstyle c_3}$}
}
\end{pspicture}
$$
A \emph{gallery\/} in a chamber system $\Delta$ is a sequence of
chambers
\begin{equation}
  \label{eq:6}
c_0\sim_{i_1} c_1\sim_{i_2}\cdots\sim_{i_k} c_k  
\end{equation}
with $c_{j-1}$ and $c_j$ $i_j$-adjacent and
$c_{j-1}\not=c_{j}$. The last condition is a technicality to help
with the accounting. We say that the gallery (\ref{eq:6}) has type
$i_1i_2\ldots i_k$, and write $c_0\rightarrow_f c_k$ where
$f=i_1i_2\ldots i_k$. If $J\subseteq I$ then a \emph{$J$-gallery\/} is
a gallery of type $i_1i_2\ldots i_k$ with the $i_j\in J$. 

A subset
$\Delta'\subseteq\Delta$ of chambers is \emph{$J$-connected\/} when
any two $c,c'\in\Delta'$ can be joined by a $J$-gallery that is
contained in $\Delta'$. The \emph{$J$-residues\/} of $\Delta$ are the
$J$-connected components and they have \emph{rank\/} $|J|$. 
Thus the chambers themselves are the rank $0$ residues. The rank $1$
residues are the equivalence classes of the equivalence relations
$\sim_i$ as $i$ runs through $I$. Call these rank $1$ residues the
\emph{panels\/} of $\Delta$. The chamber system itself has rank $|I|$.

A \emph{morphism\/} $\aa:\Delta\rightarrow\Delta'$ of chamber systems
(both over the same set $I$) is a map of the chambers of $\Delta$ to the
chambers of $\Delta'$ that preserves $i$-adjacence for all $i$: if
$c\sim_i c'$ in $\Delta$ then $\aa(c)\sim_i \aa(c')$ in $\Delta'$. An \emph{isomorphism\/} is a
bijective morphism whose inverse is also a morphism.

\begin{example}
\label{example:chambers:generic}
The local picture from \S \ref{lecture1} (below left) is a chamber system over
$I=\{1,2\}$, with chambers the edges, and two chambers
$i$-adjacent when they share a common $i$-labeled vertex. The
$\{i\}$-residues, or panels, are the pairs of edges having a $i$-labeled vertex in
common; in particular each panel contains exactly two chambers and 
there is a one-one correspondence between the
panels and the vertices:
$$
\begin{pspicture}(0,0)(\textwidth,3)
\rput(2,1.5){
\rput(0.6,0.1){\BoxedEPSF{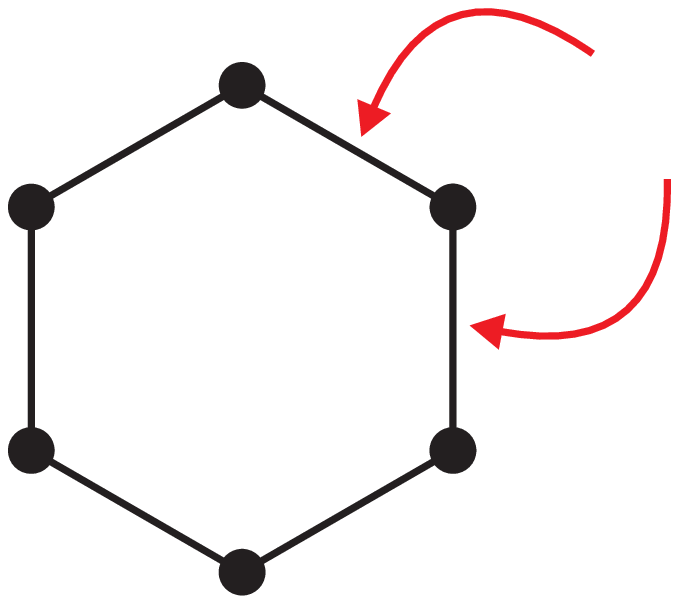 scaled 500}}
\rput(2.75,1.05){{\red $\{2\}$-residue or panel}}
\rput(-1.425,-1.525){
\rput(1.5,2.45){$\scriptstyle{1}$}
\rput(2.3,2.05){$\scriptstyle{2}$}
\rput(2.3,0.95){$\scriptstyle{1}$}
\rput(1.5,0.55){$\scriptstyle{2}$}
\rput(0.7,0.95){$\scriptstyle{1}$}
\rput(0.7,2.05){$\scriptstyle{2}$}
}
}
\rput(9,1.5){
\rput(0,0){\BoxedEPSF{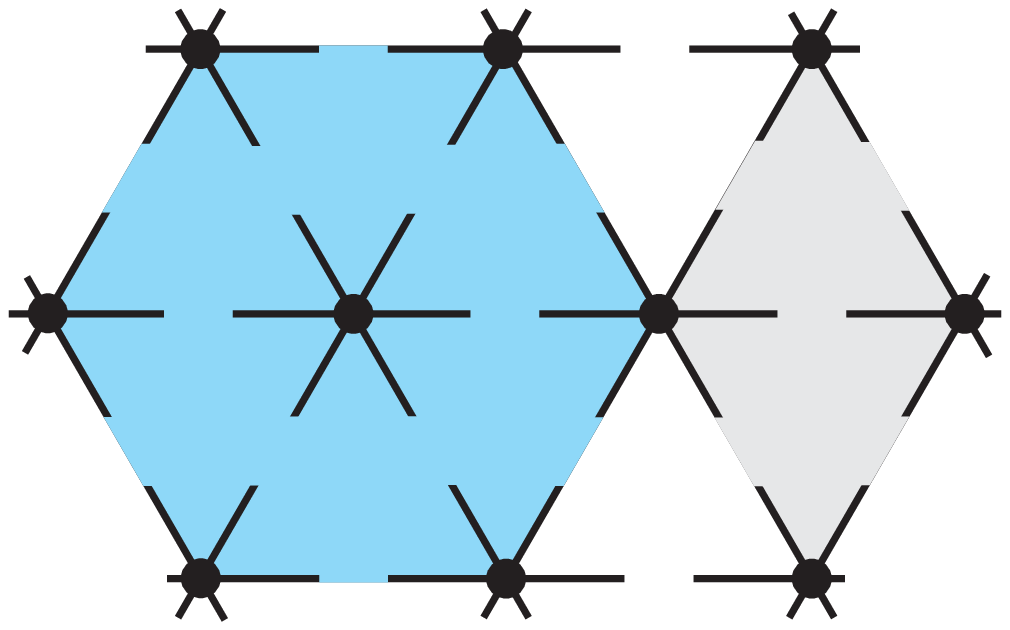 scaled 500}}
\rput(0.75,-1.35){${\scriptstyle 2}$}
\rput(0.75,1.35){${\scriptstyle 2}$}
\rput(-1.55,0){${\scriptstyle 2}$}
\rput(-0.75,-1.35){${\scriptstyle 1}$}
\rput(0,0){${\scriptstyle 3}$}
\rput(-0.75,1.35){${\scriptstyle 1}$}
\rput(1.55,0){${\scriptstyle 1}$}
\rput(0.425,0.675){${\scriptstyle 1}$}
\rput(-0.375,0.675){${\scriptstyle 2}$}
\rput(-1.2,0.675){${\scriptstyle 3}$}
\rput(-0.375,-0.675){${\scriptstyle 2}$}
\rput(-1.2,-0.675){${\scriptstyle 3}$}
\rput(-1.95,0.675){${\scriptstyle 1}$}
\rput(0.425,-0.675){${\scriptstyle 1}$}
\rput(-1.95,-0.675){${\scriptstyle 1}$}
\rput(1.95,0.675){${\scriptstyle 2}$}
\rput(1.2,0.675){${\scriptstyle 3}$}
\rput(1.95,-0.675){${\scriptstyle 2}$}
\rput(1.2,-0.675){${\scriptstyle 3}$}
}
\end{pspicture}
$$
The example above right has chambers the $2$-simplicies, $I=\{1,2,3\}$,
and two chambers $i$-adjacent when they share a 
common $i$-labeled edge. 
The six highlighted $2$-simplicies are a $\{2,3\}$-residue and the
pair of 
$2$-simplicies a $\{1\}$-residue or panel
(so again, each panel contains two chambers).
The six chambers in the rank $2$
residue have a single common vertex at their center, and there
is a one-one correspondence between the rank $2$ residues and the
vertices; similarly there is a one-one correspondence between the
panels and the edges. So the chambers are the maximal dimensional
simplicies and the residues correspond to the lower dimensional ones.
We will return to this point below.
\end{example}

\begin{example}[flag complexes]
\label{example:flagcomplexes}
Generalizing the example of \S \ref{lecture1}, let
$V$ be an $n$-dimensional vector space over a field $k$. A
\emph{flag\/} is a sequence of subspaces
$V_{i_0}\subset\cdots\subset V_{i_k}$ with $V_{i_j}$ a proper subspace
of $V_{i_{j+1}}$. Let $\Delta$ be
the chamber system over $I=\{1,\ldots,n-1\}$ whose chambers are the
\emph{maximal\/} flags $V_1\subset\cdots\subset V_{n-1}$ with $\dim V_i=i$,
and where
$$
(V_1\subset\cdots\subset V_{n-1})
\sim_i
(V'_1\subset\cdots\subset V'_{n-1})
$$
when $V_j=V'_j$ for $j\not= i$, i.e. any difference between the maximal flags occurs only in
the $i$-th position.
The chambers in the panel (or $\{i\}$-residue) containing
$V_1\subset\cdots\subset V_{n-1}$ correspond to the $1$-dimensional subspaces of the
$2$-dimensional space $V_{i+1}/V_{i-1}$. If $k$ is finite of order $q$ then each panel
thus contains $q+1$ chambers; if $k$ is infinite then each panel contains
infinitely many chambers.
\end{example}

\begin{example}[Coxeter complexes]
\label{example:coxeter.complexes}
In \S \ref{lecture2} we defined chambers, panels and $i$-adjacence for a finite
reflection group $W$ acting on a Euclidean space: the
chambers were in one-one correspondence with the elements of $W$ via
$g\leftrightarrow gc_0$ ($c_0$ a fixed fundamental
chamber), and $gc_0$ and $g'c_0$ were $i$-adjacent when $g'=gs_i$.

Now let $(W,S)$ be a Coxeter system 
with $S=\{s_i\}_{i\in I}$. The
\emph{Coxeter complex\/} $\Delta_W$ is the chamber system over $I$ with
chambers the elements of $W$ and
\begin{equation}
  \label{eq:8}
g\sim_i g'\text{ if and only if }g'=gs_i\text{ in $W$}.  
\end{equation}
Thus $g\sim_i gs_i$ and also $gs_i\sim_i gs_is_i=g$. 
The $\{i\}$-panel containing $g$ is
thus $\{g,gs_i\}$,
so each panel contains exactly two chambers (which can be thought of
as lying on either side of the panel). This is the picture the
geometry was giving us in \S \ref{lecture2}. A gallery in $\Delta_W$ has the form
$$
g\sim_{i_1} gs_{i_1}
\sim_{i_2}gs_{i_1}s_{i_2}\sim\cdots\sim_{i_k}gs_{i_1}s_{i_2}\ldots s_{i_k}.
$$
If $f=i_1i_2\ldots i_k$ 
and $s_f=s_{i_1}s_{i_2}\ldots s_{i_k}$, then there is a gallery
$g\rightarrow_f g'$ in $\Delta_W$ exactly when $g'=gs_f$ in $W$.

If $s_i,s_j\in S$ then starting at the chamber $g$ we can set
off in the two directions given by the galleries:
$$
g\sim_{i} gs_{i}\sim_{j}gs_{i}s_{j}\sim_{i}gs_{i}s_{j}s_i
\cdots
\hspace{0.5cm}
\text{ and }
\hspace{0.5cm}
g\sim_{j} gs_{j}\sim_{i}gs_{j}s_{i}\sim_{j}gs_{j}s_{i}s_j
\cdots
$$
If the order of $s_i s_j$ is finite, then $(s_i
s_j)^{m_{ij}}=1$ is equivalent to the relation
$$
s_i s_j s_i\cdots=s_j s_i s_j\cdots,
$$
where there are $m_{ij}$ symbols on both sides, so the two galleries
above, despite starting out in opposite directions, 
nevertheless end up at the same place: the chamber 
$g s_i s_j s_i\cdots=gs_j s_i s_j\cdots$. Thus the
$\{i,j\}$-residues in $\Delta_W$ are circuits containing
$2m_{ij}$ chambers when $s_is_j$ has finite order. If the order is
not finite then the residue is an infinitely long line of chambers
stretching in ``both directions'' from $g$. 
The two Coxeter groups from the end of \S \ref{lecture2} have Coxeter
complexes
illustrating both these phenomena:
$$
\begin{pspicture}(0,0)(\textwidth,3)
\rput(2.5,1.5){
\rput(0,0){\BoxedEPSF{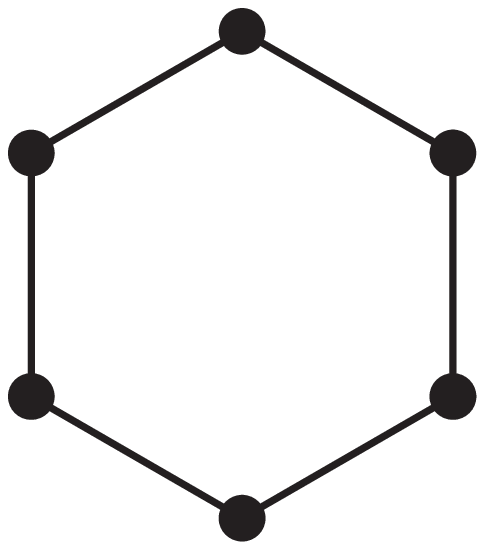 scaled 500}}
\rput(-1.475,-1.5){
\rput(1.5,2.45){$\scriptstyle{1}$}
\rput(2.3,2.05){$\scriptstyle{2}$}
\rput(2.3,0.95){$\scriptstyle{1}$}
\rput(1.5,0.55){$\scriptstyle{2}$}
\rput(0.7,0.95){$\scriptstyle{1}$}
\rput(0.7,2.05){$\scriptstyle{2}$}
}
\rput(0.6,1.175){$1$}
\rput(-0.6,1.15){$s_1$}
\rput(1.3,0){$s_2$}
\rput(-1.4,0){$s_1s_2$}
\rput(0.7,-1.15){$s_2s_1$}
\rput(-1.5,-1.15){$s_1s_2s_1=s_2s_1s_2$}
}
\rput(4,0.5){
\psline(0,0)(1,0)
\pscircle[fillstyle=solid,fillcolor=white](0,0){0.125}
\pscircle[fillstyle=solid,fillcolor=white](1,0){0.125}
\rput(0,-0.3){$s_1$}\rput(1,-0.3){$s_2$}
}
\rput(8.5,2){
\rput(0,0){\BoxedEPSF{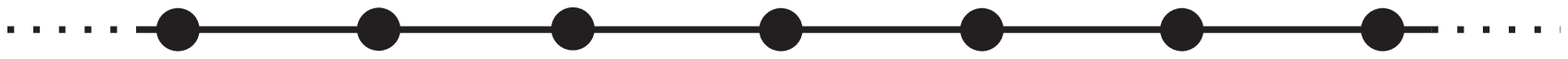 scaled 500}}
\rput(-3.3,-.3){$\scriptstyle{1}$}
\rput(-2.2,-.3){$\scriptstyle{0}$}
\rput(-1.1,-.3){$\scriptstyle{1}$}
\rput(0,-.3){$\scriptstyle{0}$}
\rput(1.1,-.3){$\scriptstyle{1}$}
\rput(2.2,-.3){$\scriptstyle{0}$}
\rput(3.3,-.3){$\scriptstyle{1}$}
\rput(-2.75,0.2){$s_0s_1s_0$}
\rput(-1.65,0.2){$s_0s_1$}
\rput(-0.555,0.2){$s_0$}
\rput(0.55,0.2){$1$}
\rput(1.65,0.2){$s_1$}
\rput(2.75,0.2){$s_1s_0$}
}
\rput(8,1){
\psline(0,0)(1,0)
\pscircle[fillstyle=solid,fillcolor=white](0,0){0.125}
\pscircle[fillstyle=solid,fillcolor=white](1,0){0.125}
\rput(0,-0.3){$s_0$}\rput(1,-0.3){$s_1$}
\rput(0.5,0.2){$\infty$}
}
\end{pspicture}
$$
\end{example}

\paragraph{Aside.}
In all our pictures of chamber systems, the chambers, panels and
lower dimensional cells have been simplicies. It turns out that chamber systems are
particularly nice examples of simplicial complexes where the chambers
are the maximal dimensional simplicies.
Moreover in all the chamber systems arising in these lectures there
is a correspondence between the lower dimensional
simplicies and the residues.

To see why, recall that an abstract simplicial
complex $X$ with vertex set $V$ is a collection of subsets of $V$
such that 
$$
\text{(a). }\ss\in X\text{ and }\tau\subset\ss\Rightarrow\tau\in X
\text{ and }
\text{(b). }\{v\}\in X\text{ for all }v\in V.
$$
A $\ss=\{v_0,\ldots,v_k\}$ is a $k$-simplex of the simplicial complex $X$.
The
empty set $\varnothing$ is by convention the unique simplex of dimension $-1$.

Now let $\Delta$ be a chamber system over $I$ and let $V$ be the set
of residues of rank $|I|-1$ (recall that there is only one residue of rank $|I|$,
namely $\Delta$ itself). Then let $X_\Delta$ be the simplicial complex
with vertex set $V$ and such that if $R_0,\ldots,R_k$ are rank
$|I|-1$ residues then
$$
\ss=\{R_0,\ldots,R_k\}\text{ is a $k$-simplex of }X_\Delta
\Leftrightarrow
\bigcap R_i\not=\varnothing.
$$
In other words, $X_\Delta$ is the \emph{nerve\/} of the covering of $\Delta$ by
rank $|I|-1$ residues.
Take the
empty intersection to be the union $\bigcup_V R_i$, and observe that the
maximum dimension a simplex can have is $|I|-1$.

If $\Delta$ is the flag complex chamber system of Example
\ref{example:flagcomplexes} with chambers the maximal flags, then the
$k$-simplicies of $X_\Delta$ correspond to the flags $V_{i_0}\subset\cdots\subset V_{i_k}$
containing $k+1$ subspaces.

We illustrate with the Coxeter complex $\Delta_W$ of the Coxeter system $(W,S)$
with the symbol shown below left.
Some elements
of $W$ have been written down in a suggestive pattern, grouped into
three rank $2$ residues. The simplicial complex $X_\Delta$ acquires a $2$-simplex
from these residues as any two intersect in a residue of rank
$1$ and all three intersect in a residue of rank $0$. In fact $X_W$ is the
infinite tiling of the plane from Example \ref{example:chambers:generic}:
$$
\begin{pspicture}(0,0)(\textwidth,5)
\rput(1,1){
\rput(0,0){\BoxedEPSF{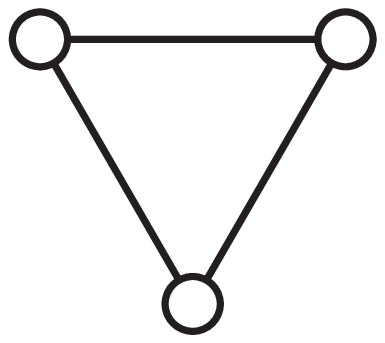 scaled 450}}
\rput(-0.7,0.95){$s_1$}\rput(0,-0.95){$s_2$}\rput(0.7,0.95){$s_3$}
\rput(0,1.3){$(W,S)$}
}
\rput(5,2.5){
\rput(0,0){\BoxedEPSF{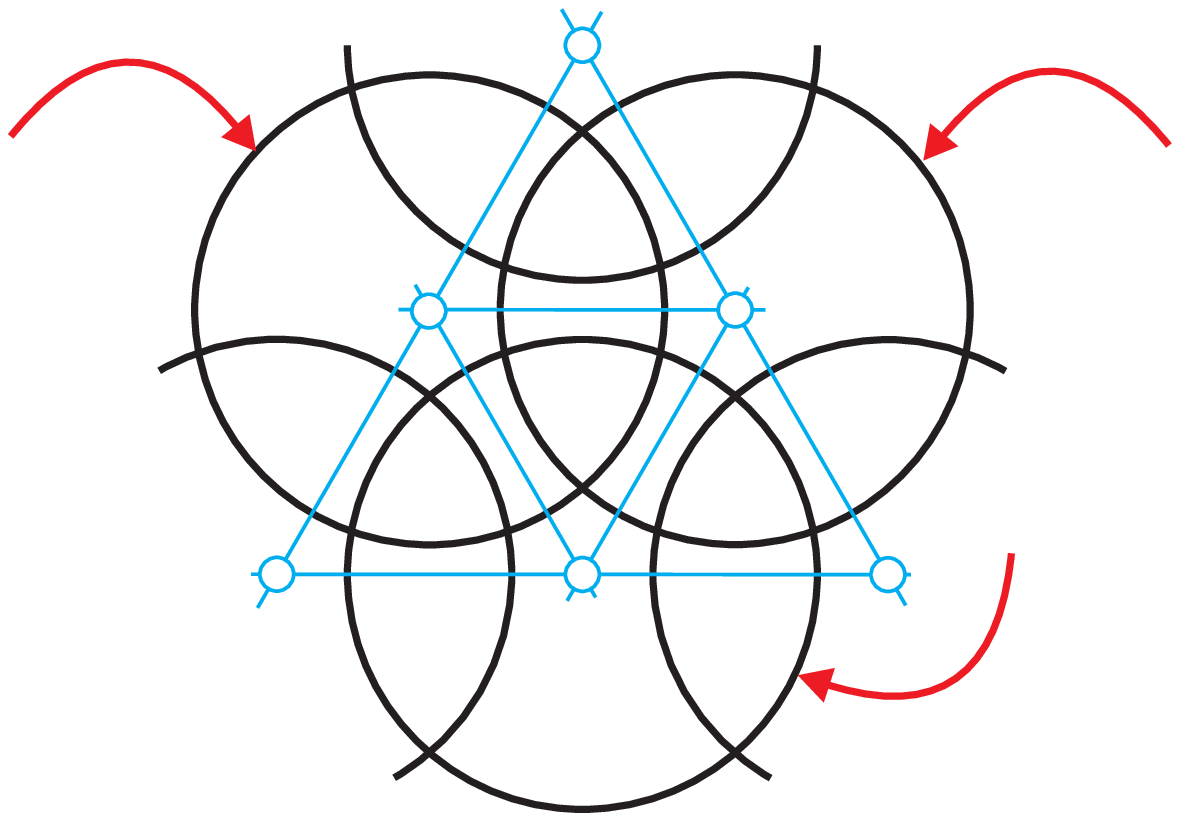 scaled 650}}
\rput(0,1.25){$s_1s_3$}
\rput(-1,-0.55){$s_1s_2$}
\rput(1,-0.55){$1$}
\rput(0,0){$s_1$}
\rput(2,0){$s_3$}
\rput(2,1.25){$s_3s_1$}
\rput*(1,1.85){$s_1s_3s_1$}
\rput*(-1,1.85){$s_1s_3s_2$}
\rput*(-2,0){$s_1s_2s_3$}
\rput*(-2,1.25){$s_1s_2s_3s_2$}
\rput(1,-1.7){$s_2$}
\rput*(-1,-1.7){$s_1s_2s_1$}
\rput(0,-2.3){$s_2s_1$}
\rput(4.2,1.5){{\red $\{1,3\}$-residue}}
\rput(-4,1.5){{\red $\{2,3\}$-residue}}
\rput(3.5,-0.7){{\red $\{1,2\}$-residue}}
\rput(-2.25,-1.75){$\Delta_W$}
}
\rput(11,2){
\rput(0,0){\BoxedEPSF{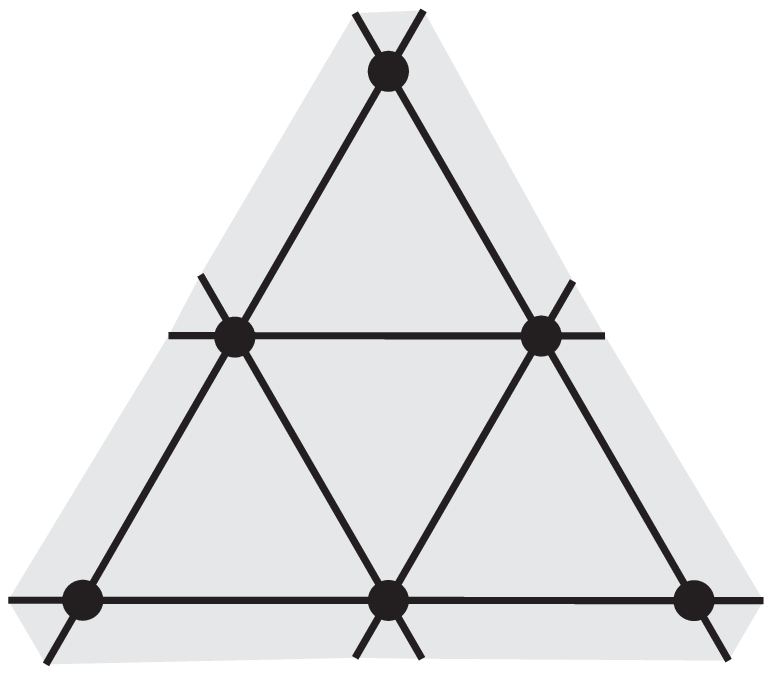 scaled 500}}
\rput(1,1.5){$X_\Delta$}
}
\end{pspicture}
$$
It would seem from this example that if
$R_0,\ldots,R_k$ are rank $|I|-1$ residues over $J_0,\ldots,J_k$
with $\bigcap R_i\not=\varnothing$, then $\bigcap R_i$ is a residue
over $\bigcap J_i$. In fact this is always true for a Coxeter complex and
indeed any building, although not for an arbitrary chamber system.
As $\bigcap J_i$ has $|I|-(k+1)$ elements, there is a
one-one correspondence between the simplicies of $X_\Delta$ and the
residues of $\Delta$:
$$
\text{codimension }\ell\text{ simplicies }\ss=\{R_0,\ldots,R_{m}\}
\leftrightarrow
\text{ residues }\bigcap_{i=0}^{m} R_i\text{ of rank }\ell, 
$$
where $m=|I|-(\ell+1)$. So for buildings the chambers of a chamber
system $\Delta$ are the top
dimensional simplicies of $X_\Delta$, with the lower dimensional simplicies
given by the residues.

\paragraph{}
Returning to the general discussion, we now have all the properties of chamber systems that we
need. We finish the section by defining a $W$-valued metric
on a Coxeter complex $\Delta_W$. 

If $(W,S)$ is a Coxeter system and $f=i_1i_2\ldots i_k$ 
with $s_f=s_{i_1}s_{i_2}\ldots s_{i_k}$, then we have seen that there is a gallery
$g\rightarrow_f g'$ in $\Delta_W$ exactly when $g'=gs_f$ in $W$. Call
such a gallery \emph{minimal\/} if there is no gallery in $\Delta_W$ from $g$ to
$g'$ that passes through fewer chambers. Call an expression
$s_f=s_{i_1}s_{i_2}\ldots s_{i_k}$ \emph{reduced\/} if there is no
expression in $W$ for $s_f$ involving fewer $s$'s (counted with
multiplicity). Thus a gallery $g\rightarrow_f g'$ is minimal if and
only if the expression $s_f$ is reduced. 

Define $\delta_W:\Delta_W\times\Delta_W\rightarrow W$ by
$\delta_W(g,g')=g^{-1}g'$. Then
\begin{equation}
  \label{eq:7}
\delta_W(g,g')=s_f
\Leftrightarrow  
g'=gs_f
\Leftrightarrow
\text{ there is a gallery }g\rightarrow_f g'.
\end{equation}
Moreover, $\delta_W(g,g')$ is reduced if and only if the gallery 
$g\rightarrow_f g'$ is minimal. A slight relaxation will define 
the metric on an arbitrary building.
Here are two examples, one of which is our running one:
$$
\begin{pspicture}(0,0)(\textwidth,4)
\rput(3.8,2.25){
\rput(0,0){\BoxedEPSF{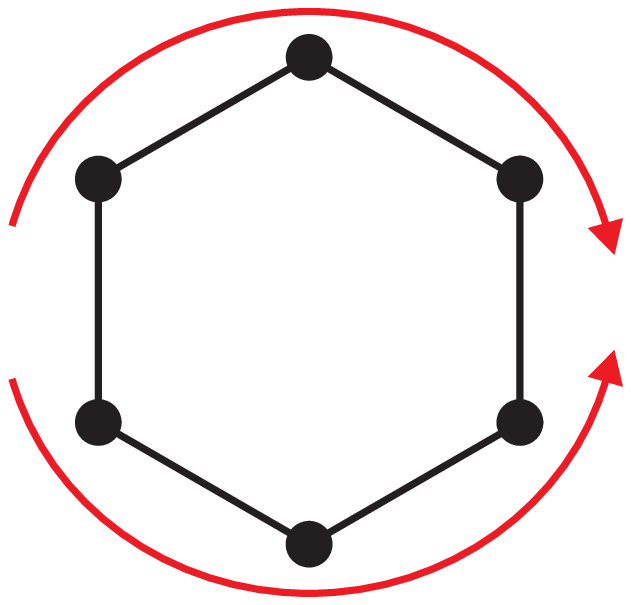 scaled 500}}
\rput(-1.475,-1.5){
\rput(1.5,2.45){$\scriptstyle{1}$}
\rput(2.3,2.05){$\scriptstyle{2}$}
\rput(2.3,0.95){$\scriptstyle{1}$}
\rput(1.5,0.55){$\scriptstyle{2}$}
\rput(0.7,0.95){$\scriptstyle{1}$}
\rput(0.7,2.05){$\scriptstyle{2}$}
}
\rput(1.75,0){$g'=s_2$}
\rput(-1.9,0){$g=s_1s_2$}
}
\rput(0.5,3.5){
\psline(0,0)(1,0)
\pscircle[fillstyle=solid,fillcolor=white](0,0){0.125}
\pscircle[fillstyle=solid,fillcolor=white](1,0){0.125}
\rput(0,-0.3){$s_1$}\rput(1,-0.3){$s_2$}
}
\rput(2.2,0.5){$\delta_W(g,g')=s_2s_1s_2\,\,(=s_1s_2s_1 )$}
\rput(10.5,2.25){
\rput(0,0){\BoxedEPSF{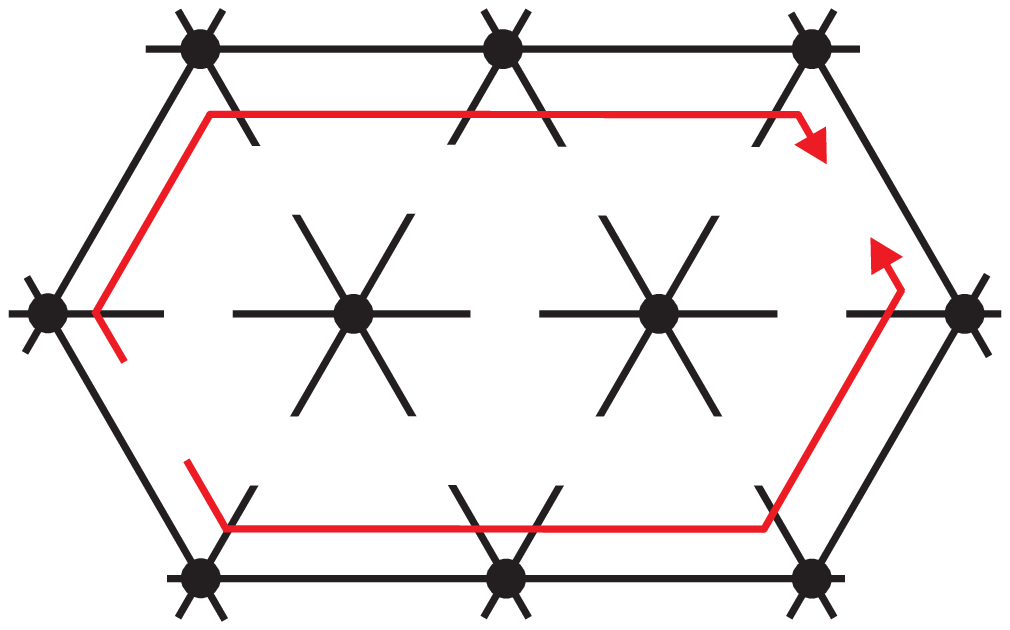 scaled 500}}
\rput(-1.55,0){${\scriptstyle 2}$}
\rput(0,0){${\scriptstyle 3}$}
\rput(1.55,0){${\scriptstyle 1}$}
\rput(0.425,0.675){${\scriptstyle 1}$}
\rput(-0.375,0.675){${\scriptstyle 2}$}
\rput(-1.2,0.675){${\scriptstyle 3}$}
\rput(-0.375,-0.675){${\scriptstyle 2}$}
\rput(-1.2,-0.675){${\scriptstyle 3}$}
\rput(0.425,-0.675){${\scriptstyle 1}$}
\rput(1.2,0.675){${\scriptstyle 3}$}
\rput(1.2,-0.675){${\scriptstyle 3}$}
\rput(-1.6,-0.5){$g$}\rput(1.45,0.4){$g'$}
}
\rput(9.8,0.5){$\delta_W(g,g')=s_2s_3s_2s_1s_3\,\,(=s_3s_2s_1s_3s_1=\text{etc} )$}
\rput(7.25,3){
\rput(0,0){\BoxedEPSF{fig22.eps scaled 450}}
\rput(-0.7,0.95){$s_1$}\rput(0,-0.95){$s_2$}\rput(0.7,0.95){$s_3$}
}
\end{pspicture}
$$

\begin{figure}
  \centering
\begin{pspicture}(0,0)(12,6.5)
\rput(2.5,3.25){\BoxedEPSF{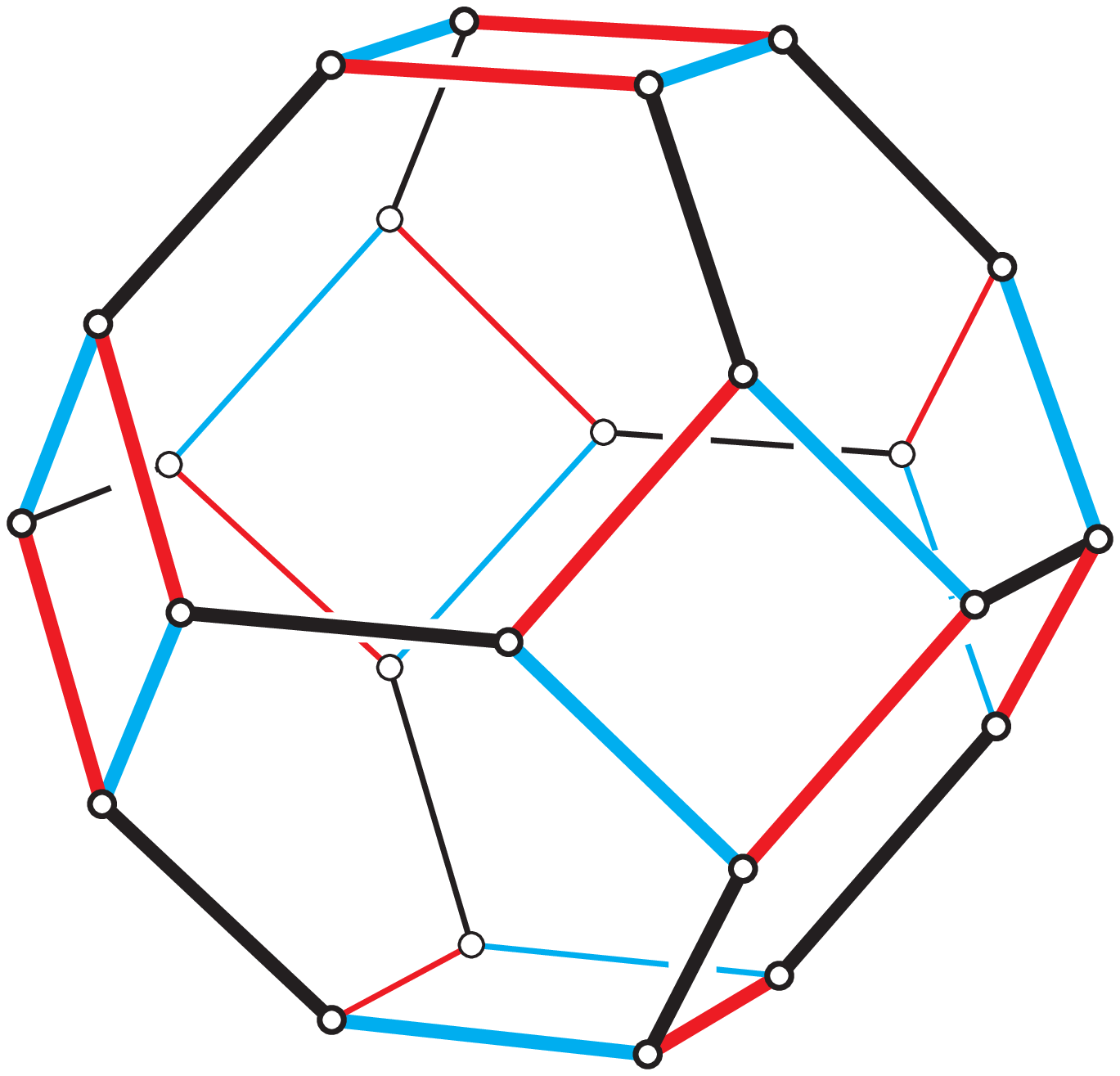 scaled 350}}
\rput(9,3.25){\BoxedEPSF{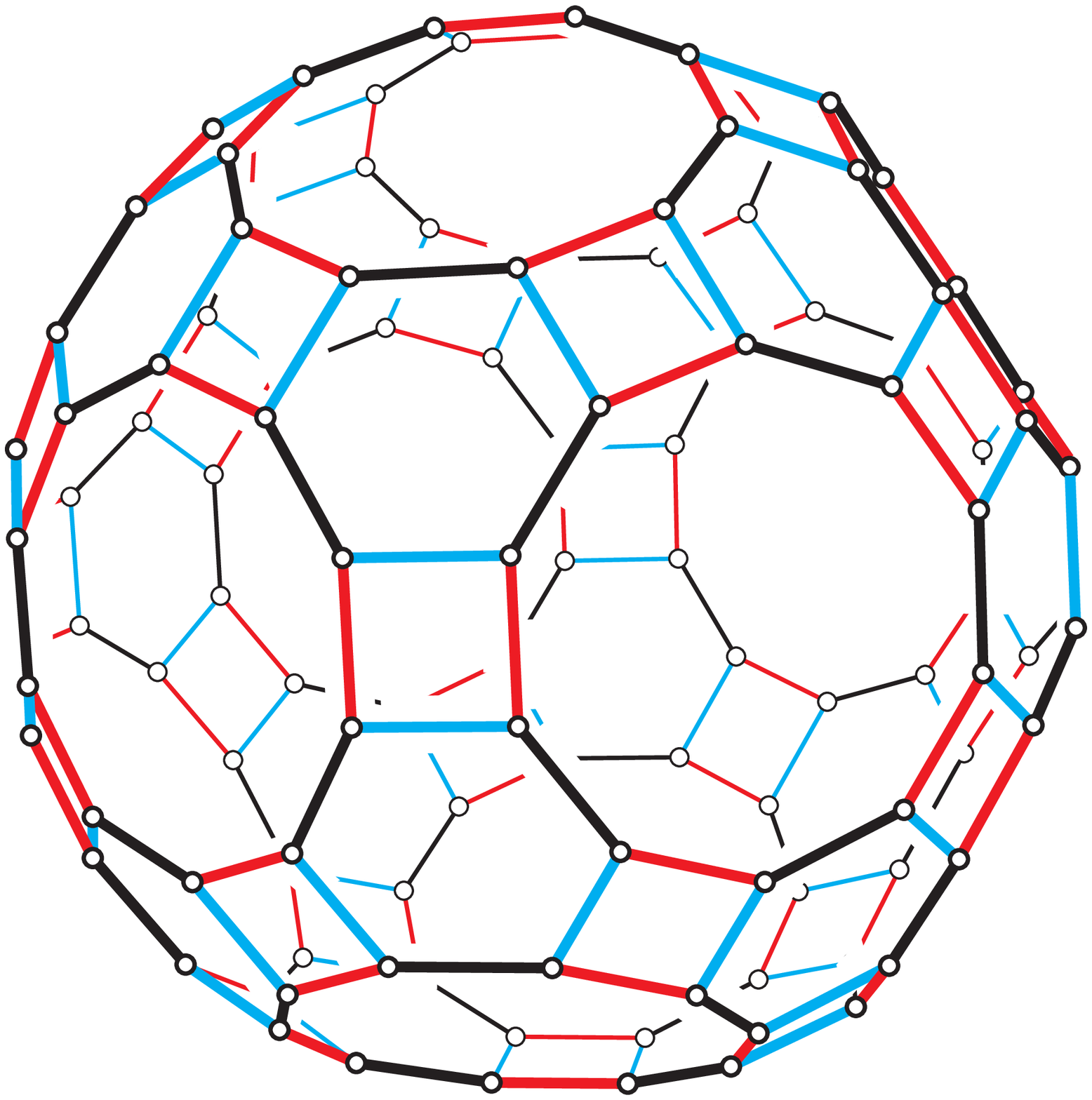 scaled 350}}
\end{pspicture}
\caption{Chamber systems as edge coloured graphs. The local picture
  for the flag complex of a four dimensional space \emph{(left)} and the
  Coxeter complex of the group of symmetries of a dodecahedron
  \emph{(right)}. Both are chamber systems over $I=\{1,2,3\}$.}
  \label{fig:permutohedron}
\end{figure}

\paragraph{Another way to draw chamber systems.} A chamber system over
$I$ can be drawn as a graph whose edges are ``coloured'' by $I$. The
vertices of the graph are the chambers, and two vertices are joined by an
edge labeled $i\in I$ 
iff the corresponding chambers are $i$-adjacent. These graphs are
essentially the $1$-skeletons of the duals of our simplicial
complexes. If $\Delta_W$ is the Coxeter complex of the Coxeter system
$(W,S)$ then this graph is the Cayley graph of $W$ with respect to the
generating set $S$. Figure \ref{fig:permutohedron} (left) shows the graph for
the local picture of the flag complex of a four dimensional space of
Figure \ref{fig:localpicture2} (or the Coxeter 
complex of \begin{pspicture}(0,0)(2.5,0.3)
\rput(0.25,0.1){
\psline(0,0)(1,0) \psline(1,0)(2,0)
\pscircle[fillstyle=solid,fillcolor=white](0,0){0.125}
\pscircle[fillstyle=solid,fillcolor=white](1,0){0.125}
\pscircle[fillstyle=solid,fillcolor=white](2,0){0.125}
}
\end{pspicture}) and (right) the graph for the Coxeter complex of
the group of symmetries of the dodecahedron (with Coxeter symbol \begin{pspicture}(0,0)(2.5,0.3)
\rput(0.25,0.1){
\psline(0,0)(1,0) \psline(1,0)(2,0)
\pscircle[fillstyle=solid,fillcolor=white](0,0){0.125}
\pscircle[fillstyle=solid,fillcolor=white](1,0){0.125}
\pscircle[fillstyle=solid,fillcolor=white](2,0){0.125}
\rput(1.5,0.15){${\scriptstyle 5}$}
}
\end{pspicture}).

\section{Buildings and Apartments}
\label{lecture4}

Let $(W,S)$ be a Coxeter system with $S=\{s_i\}_{i\in I}$. 
A \emph{building of type $(W,S)$\/} is a chamber system $\Delta$ over
$I$ such that:
\begin{description}
\item[(B1).] every panel of $\Delta$ contains at least two chambers;
\item[(B2).] $\Delta$ has a $W$-valued metric
  $\delta:\Delta\times\Delta\rightarrow W$ such that if
  $s_f=s_{i_1}\ldots s_{i_k}$ is a \emph{reduced\/} expression in $W$ then 
$$
\delta(c,c')=s_f\Leftrightarrow\mbox{ there is a gallery
}c\rightarrow_f c'\mbox{ in }\Delta.
$$
\end{description}

\begin{example}[Coxeter complexes]
There is at least one building for every Coxeter system $(W,S)$, namely
the Coxeter complex $\Delta_W$ with $\delta=\delta_W$ in (\ref{eq:7}),
hence (B2). For (B1) we observed in Example \ref{example:coxeter.complexes}
that
the panels in $\Delta_W$ have the form $\{g,gs\}$ for $g\in W$ and $s\in S$. Such a
building, where each panel has the minimum possible number of
chambers, is said to be \emph{thin\/}. It turns out that the thin
buildings are precisely the Coxeter complexes.
\end{example}

\begin{example}[a spherical building of type \begin{pspicture}(0,0)(4.1,0.4)
\rput(0.25,0.1){
\psline(0,0)(1,0)
\psline(1,0)(1.4,0)
\psline[linestyle=dotted](1.5,0)(2.2,0)
\psline(2.3,0)(2.7,0)
\psline(2.7,0)(3.7,0)
\pscircle[fillstyle=solid,fillcolor=white](0,0){0.125}
\pscircle[fillstyle=solid,fillcolor=white](1,0){0.125}
\pscircle[fillstyle=solid,fillcolor=white](2.7,0){0.125}
\pscircle[fillstyle=solid,fillcolor=white](3.7,0){0.125}
}
\end{pspicture}]
\label{example:Anbuilding}
For $(W,S)$ having this symbol ($n-1$ vertices)
we put a $W$-valued metric on the flag complex of Example
\ref{example:flagcomplexes}. First identify $(W,S)$ with $\gS_n$ as in
\S\ref{lecture2}, with $s_i\mapsto (i,i+1)$ for $1\leq i\leq n-1$. Let
$$
c=(V_1\subset\cdots\subset V_{n-1})
\text{ and }
c'=(V'_1\subset\cdots\subset V'_{n-1})
$$
be chambers and write $V_0=V'_0=0$, $V_n=V'_n=V$. We can define
$\delta(c,c')\in\gS_n$ using the filtration of
$V'_i/V'_{i-1}$ of \S\ref{lecture1} in the obvious way. Alternatively,
for $1\leq i\leq n$, let
$$
\pi(i)=\min\{j\,|\,V'_i\subset V'_{i-1}+V_j\}
$$
and define $\delta(c,c')=\pi$. We show that we have a building (when
$\dim V=3$) at the end of this section.
\end{example}

\begin{example}[an affine building of type \begin{pspicture}(0,0)(1.5,0.4)
\rput(0.25,0.1){
\psline(0,0)(1,0)
\pscircle[fillstyle=solid,fillcolor=white](0,0){0.125}
\pscircle[fillstyle=solid,fillcolor=white](1,0){0.125}
\rput(0.5,0.15){$\infty$}
}
\end{pspicture}]
\label{example:affinebuilding}
An affine building has type $(W,S)$ an
affine reflection group as in Example \ref{example:affine}.
Taking this example, 
with $S=\{s_0,s_1\}$ and Coxeter symbol 
\begin{pspicture}(0,0)(1.5,0.4)
\rput(0.25,0.1){
\psline(0,0)(1,0)
\pscircle[fillstyle=solid,fillcolor=white](0,0){0.125}
\pscircle[fillstyle=solid,fillcolor=white](1,0){0.125}
\rput(0.5,0.15){$\infty$}
}
\end{pspicture},
let $\Delta$ be the
chamber system over $I=\{0,1\}$ shown below -- an infinite 3-valent tree. 
The edges are the chambers, and two chambers
are $0$-adjacent when they share a common black vertex and
$1$-adjacent when they share a common white vertex. 
Each panel thus contains three chambers, hence (B1).
The Coxeter complex $\Delta_W$ is in Example \ref{example:coxeter.complexes} (also
a tree).
$$
\begin{pspicture}(0,0)(\textwidth,5)
\rput(3.25,0){
\rput(3,2.5){
\rput(0,0){\BoxedEPSF{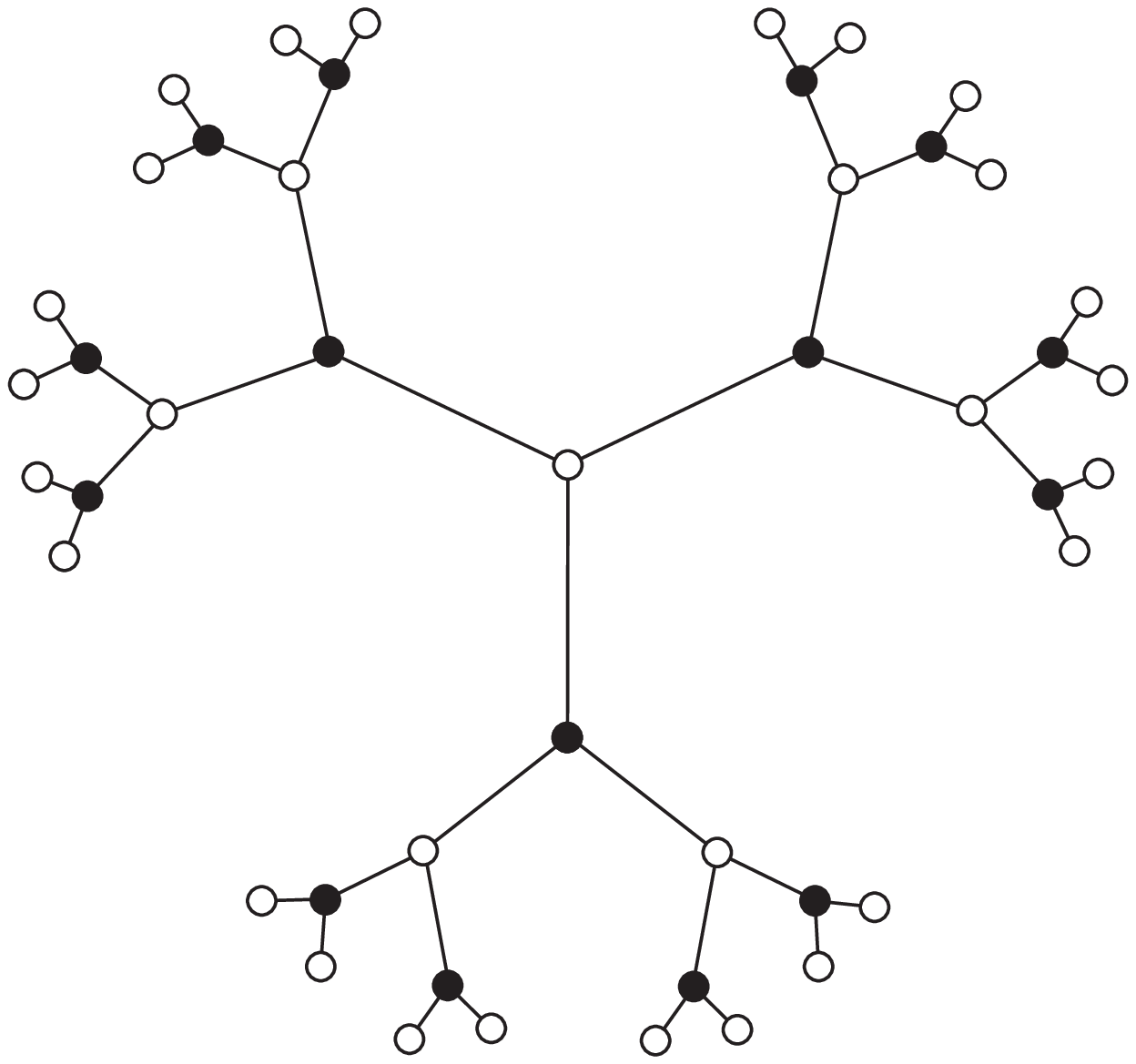 scaled 450}}
}
\rput(6.5,2.75){
\rput(0,0){\BoxedEPSF{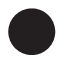 scaled 450}}
}
\rput(6.5,2.25){
\rput(0,0){\BoxedEPSF{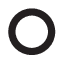 scaled 450}}
}
\rput(7.7,2.8){$=0$-adjacent}\rput(7.7,2.3){$=1$-adjacent}
\rput(3,4){$\Delta$}
}
\end{pspicture}
$$
To define the $W$-metric on $\Delta$ recall that in a tree there is a
unique path between chambers without ``backtracking'': a backtrack is a
path that crosses an edge and then immediately comes back across the
edge again. For chambers $c,c'\in\Delta$,
match this unique path between $c$ and $c'$ with the same path
starting at $1$ in the Coxeter complex $\Delta_W$:
$$
\begin{pspicture}(0,0)(\textwidth,3)
\rput(6.25,1.5){
\rput(0,0){\BoxedEPSF{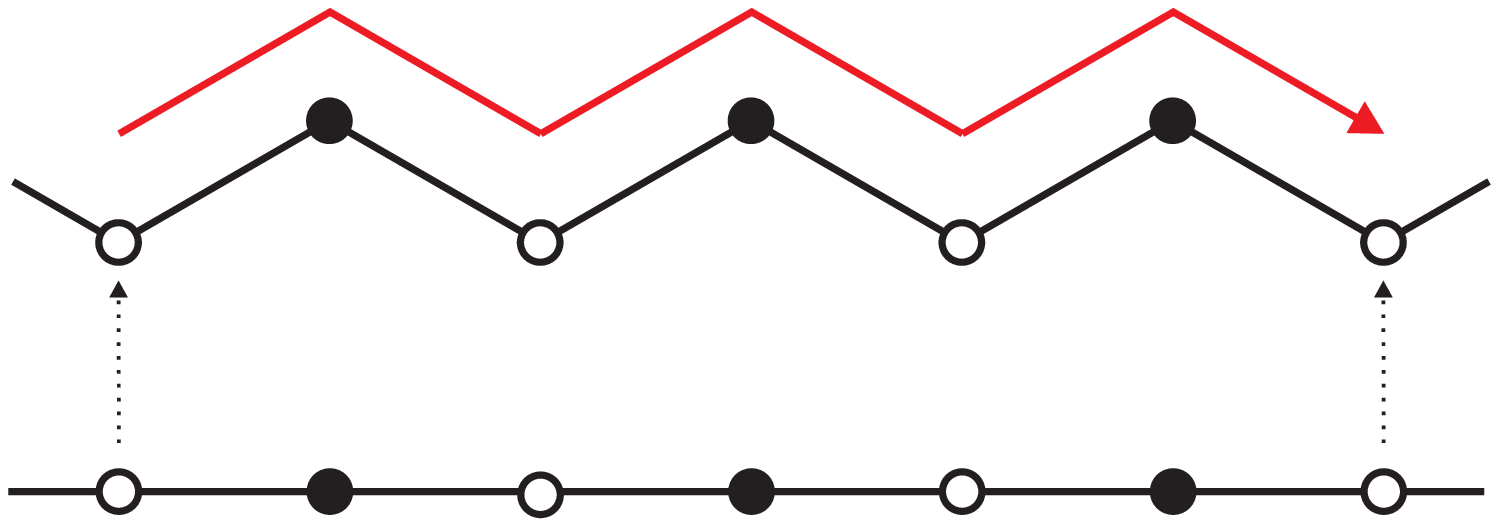 scaled 500}}
}
\rput(6.25,3){{\red unique path}}
\rput(11.25,2){$\Delta$}\rput(11.25,0.3){$\Delta_W$}
\rput(3.6,0.6){$1$}\rput(9,0.6){$g$}
\rput(3.6,1.7){$c$}\rput(8.9,1.7){$c'$}
\end{pspicture}
$$
and define $\delta(c,c')$ to be the resulting $g$. 
To see (B2), let $\delta(c,c')=g\in W$ and suppose that
$g=s_{j_1}\ldots s_{j_\ell}$. Then by 
(\ref{eq:7}) there is a gallery in $\Delta_W$
from $1$ to $g$ of type $j_1\ldots j_\ell$. As $\Delta_W$ is also a
tree this gallery differs from the unique minimal one only by
backtracks. First transfer this minimal gallery to $\Delta$ to get 
the minimal gallery from $c$ to $c'$, and then transfer the backtracks
to obtain a gallery of type $j_1\ldots j_\ell$ from $c$ to
$c'$. Conversely if there is a gallery from $c$ to $c'$ of type
$j_i\ldots j_\ell$ with $s_{j_1}\ldots s_{j_\ell}$ \emph{reduced\/}, then in
particular no two consecutive $s$'s are the same and so the gallery has no
backtracks. Thus it is \emph{the\/} unique minimal gallery from $c$ to
$c'$ giving $\delta(c,c')=s_{j_1}\ldots s_{j_\ell}$ by definition.
\end{example}

In a Coxeter complex we have $\delta_W(c,c')=s_{i_1}\ldots s_{i_k}$ if and
only if there is a gallery of type $i_1\ldots i_k$ from $c$ to $c'$,
but in an arbitrary building there is the extra condition that the
word $s_{i_1}\ldots s_{i_k}$ be reduced. We can see why in the example
above: if there is a gallery of type $i_1\ldots i_k$ from $c$ to $c'$
with $s_{i_1}\ldots s_{i_k}$ not reduced, then $\delta(c,c')$ need not
necessarily be $s_{i_1}\ldots s_{i_k}$. For example, if we have three
adjacent chambers:
$$
\begin{pspicture}(0,0)(3,2.75)
\rput(0,-.25){
\rput(1.5,1.5){
\rput(0,0){\BoxedEPSF{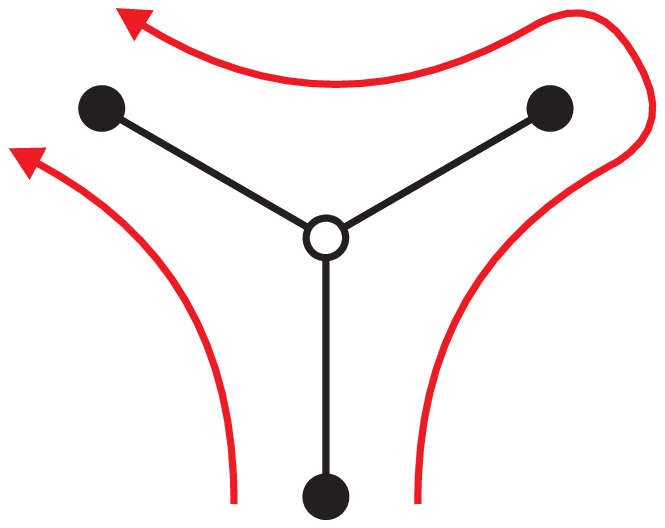 scaled 500}}
}
\rput(1.3,1){$c$}\rput(1,2.2){$c'$}
\rput(0.4,1){${\scriptstyle{\red f=1}}$}
\rput(2.5,1){${\scriptstyle{\red f'=11}}$}}
\end{pspicture}
$$
then there is a gallery of type $1$ from
$c$ to $c'$ with $s_1$ reduced, hence $\delta(c,c')=s_1$. The
non-reduced gallery $c\rightarrow_{11} c'$ does not give $\delta(c,c')=s_1s_1$,
as $s_1s_1=1\not=s_1$. 

Examples \ref{example:Anbuilding}-\ref{example:affinebuilding} are our first of  \emph{thick\/}
buildings: one where every panel contains at least three
chambers. ``Thick'' is generally taken to be synonymous with
interesting.

It turns out that there are quite naturally arising Coxeter groups
for which there are \emph{no\/} thick buildings. One such example is
the group of reflectional symmetries of a regular dodecahedron having
symbol \begin{pspicture}(0,0)(2.5,0.3)
\rput(0.25,0.1){
\psline(0,0)(1,0) \psline(1,0)(2,0)
\pscircle[fillstyle=solid,fillcolor=white](0,0){0.125}
\pscircle[fillstyle=solid,fillcolor=white](1,0){0.125}
\pscircle[fillstyle=solid,fillcolor=white](2,0){0.125}
\rput(1.5,0.15){${\scriptstyle 5}$}
}
\end{pspicture}.

In \S \ref{lecture1} (as well as Example 
\ref{example:affinebuilding}) we defined the $W$-metric $\delta$ by situating a pair of
chambers $c,c'$ inside a copy of the Coxeter complex $\Delta_W$ and
transferring the metric $\delta_W$ defined in (\ref{eq:7}). 
We need to see that this process is well defined -- although this is
obvious in Example \ref{example:affinebuilding} -- and that the
resulting $\delta$ satisfies (B2). This leads to an
alternative definition of building (Theorem \ref{theorem:second_building_definition}
below) based on this idea of defining
$\delta$ locally.

Let $(\Delta,\delta)$ and $(\Delta',\delta')$ be buildings of type
$(W,S)$ and $X\subset (\Delta,\delta),Y\subset (\Delta',\delta')$ be
subsets. A morphism $\aa:X\rightarrow Y$ is an \emph{isometry\/} when it
preserves the $W$-metrics: for all chambers $c,c'$ in $X$ we have 
$\delta'(\aa(c),\aa(c'))=\delta(c,c')$. A simple example is if $g_0\in
W$, then $g\mapsto g_0g$ is an isometry 
$\Delta_W\rightarrow\Delta_W$.

The following result guarantees the existence of copies of the Coxeter
complex in a building:

\begin{theorem}
\label{theorem:isometries}
Let $\Delta$ be a building of type $(W,S)$ and $X$ a subset of the
Coxeter complex $\Delta_W$. Then any isometry $X\rightarrow\Delta$
extends to an isometry $\Delta_W\rightarrow\Delta$. 
\end{theorem}

An \emph{apartment\/} in a building $\Delta$ of type $(W,S)$ is an isometric image of the
Coxeter complex $\Delta_W$, i.e. a subset of the form $\aa(\Delta_W)$ for 
$\aa:\Delta_W\rightarrow\Delta$ some isometry.
Apartments are precisely the local pictures we saw in \S \ref{lecture1}.

We are particularly interested in the following two consequences of
Theorem \ref{theorem:isometries}:
\begin{equation}
  \label{eq:9}
\text{Any two chambers $c,c'$ lie in some apartment $A$.}  
\end{equation}
(If $\delta(c,c')=g\in W$,
then $X=(1,g)\subset\Delta_W\mapsto(c,c')\subset\Delta$ is an
isometry. It extends by Theorem \ref{theorem:isometries} to an
isometry $\Delta_W\rightarrow\Delta$ and hence an apartment containing
$c,c'$.) 
So the $W$-metric on $\Delta$ can be recovered from the metric on the
Coxeter complex; moreover, the metrics on overlapping Coxeter complexes
agree on the overlaps:
\begin{equation}
  \label{eq:10}
\text{If chambers $c,c'\in A$ and $c,c'\in A'$
 then there is an isometry $A\rightarrow A'$ 
 fixing $A\cap A'$.}  
\end{equation}
(We leave this to the reader with the following hints: use
the apartments to get an isometry $A\rightarrow A'$ fixing a
chamber $c_0\in A\cap A'$; then show that every chamber in the
intersection is fixed by showing that in an apartment there is a
unique chamber a given $W$-distance from $c_0$.)


It turns out that any chamber system covered by sufficiently many
Coxeter complexes in a sufficiently nice way so that (\ref{eq:9}) and
(\ref{eq:10}) hold can be made into a building by patching together
the local metrics on the Coxeter complexes \emph{ala\/} \S
\ref{lecture1}. 

To
formulate this properly we need to replace isometries by maps not
involving metrics.
Let $\Delta,\Delta'$ be chamber systems over the
same set $I$. We leave it as an exercise to show that (i). 
$\aa:(\Delta,\delta)\rightarrow (\Delta',\delta')$ is an isometry of
buildings if
and only if $\aa:\Delta\rightarrow\Delta'$ is an injective morphism of
chamber systems, and (ii). $\aa$ is a surjective isometry of buildings if
and only if $\aa$ an isomorphism of chamber systems. 

\begin{theorem}
\label{theorem:second_building_definition}
Let $(W,S)$ be a Coxeter system with $S=\{s_i\}_{i\in I}$ and
$\Delta$ a chamber system over $I$. Suppose $\Delta$ contains a collection $\{A_\aa\}$ of
sub-chamber systems over $I$, called
apartments, with each subsystem isomorphic (as a chamber system) to the
Coxeter complex $\Delta_W$.
Suppose also that
\begin{description}
  \item[(B1$^\prime$).] any two chambers $c,c'$ of $\Delta$ are contained in some
    apartment $A$, and
  \item[(B2$^\prime$).] if chambers $c,c'\in A_\aa$ and $\in A_\bb$,
    then there is an isomorphism $A_\aa\rightarrow A_\bb$ fixing $A_\aa\cap A_\bb$. 
\end{description}
Define $\delta:\Delta\times\Delta\rightarrow W$ by
$\delta(c,c')=\delta_W(\aa(c),\aa(c'))$ where $\aa:\Delta_W\rightarrow
A$ is an isomorphism with $c,c'\in A$. Then $(\Delta,\delta)$ is a
building of type $(W,S)$. 
\end{theorem}

\begin{example}[the flag complex of \S \ref{lecture1} revisited]
The chamber system structure on the flag complex $\Delta$ of \S \ref{lecture1}
was given there (and in Example \ref{example:flagcomplexes}, where we
saw that $\Delta$ is thick). 
If $L_1,L_2,L_3$ are lines in $V$ spanned by independent vectors, then
we get a hexagonal configuration as in \S \ref{lecture1}. Let the apartments
be all the hexagons obtained in this way. If 
$c=V_1\subset V_2$ and $c'=V_1'\subset V_2'$ are chambers, then they
can be situated in an apartment by extending $V_1,V_1'$ to a set
$L_1,L_2,L_3$ of independent lines. If 
$V_1\not=V_1'$, $V_2\not=V_2'$ and
$V_2\cap V_2'$ is a line different from $V_1,V_1'$ as for the $c,c'$
of \S \ref{lecture1}, then this extension is unique, so $c,c'$ lie in a unique
apartment. Otherwise (e.g. if $V_2\cap V_2'$ is one of $V_2$ or
$V_2'$) there is some choice. In any case, if $L_1,L_2,L_3$  and
$L'_1,L'_2,L'_3$ are two such extensions corresponding to apartments
$A_\aa,A_\beta$ containing $c,c'$, then any $g\in GL(V)$ with $g(L_i)=L'_i$ induces an
isomorphism $A_\aa\rightarrow A_\beta$ that fixes $A_\aa\cap A_\beta$.
\end{example}

\section{Spherical Buildings}
\label{lecture5}


So far our supply of \emph{thick\/}
buildings is a little disappointing: only 
the flag complex of \S \ref{lecture1} and the affine building of
Example \ref{example:affinebuilding}. In this section we considerably
increase the library by extracting a building from the
structure of a reductive algebraic group. These guys really are the motivating
examples of buildings.

Call a building of type $(W,S)$ \emph{spherical\/} when the Coxeter system
$(W,S)$ is spherical (i.e. finite). 
It turns out that there is a uniform construction of a large class of
thick spherical buildings.
To motivate this we reconstruct the flag complex
building $\Delta$ of \S \ref{lecture1} inside the general linear group
$G=GL(V)\cong GL_3(k)$.

First, let $P\subset G$ be the subgroup of
permutation matrices -- those matrices with exactly one $1$ in each row and
column and all other entries $0$; alternatively, the
$a_\pi=\sum_j e_{\pi j,j}$, where $\pi\in\gS_3$ and $e_{ij}$ is the $3\times 3$
matrix with a $1$ in the $ij$-th position and $0$'s elsewhere. 
The map $\pi\mapsto a_\pi$ is an
isomorphism $\gS_3\rightarrow P$ with 
\begin{equation}
  \label{eq:12}
s_1=(1,2)\mapsto
\left(\begin{array}{ccc}
  0&1&0\\
  1&0&0\\
  0&0&1
\end{array}\right)
\text{ and }
s_2=(2,3)\mapsto
\left(\begin{array}{ccc}
  1&0&0\\
  0&0&1\\
  0&1&0
\end{array}\right).  
\end{equation}
For the rest of this section we will blur the distinction between the
symmetric group $\gS_3$, the group of permutation matrices $P$, and the
Coxeter system $(W,S)$ with the symbol 
\begin{pspicture}(0,0)(1.5,0.3)
\rput(0.25,0.1){
\psline(0,0)(1,0)
\pscircle[fillstyle=solid,fillcolor=white](0,0){0.125}
\pscircle[fillstyle=solid,fillcolor=white](1,0){0.125}
}
\end{pspicture}.

Assume for the moment that:
\begin{description}
  \item[(G1).] The action of $G$ on the flag complex $\Delta$
    given by $a:V_1\subset V_2\mapsto aV_1\subset aV_2$ for $a\in G$, is by chamber
    system isomorphisms (hence via isometries by the comments
    immediately prior to Theorem \ref{theorem:second_building_definition}).

\item[(G2).] Fix $g\in (W,S)$ and let
  $X(g)=\{(c,c')\in\Delta\times\Delta\,|\,\delta(c,c')=g\}$. Then for
  any $g$
  the diagonal action $a:(c,c')\mapsto (ac,ac')$ of $G$ on $X(g)$ is
  transitive (thus $G$ acts transitively on the ordered pairs
  of chambers a fixed $W$-distance apart).

\item[(G3).] Let $A_0\subset\Delta$ be the apartment given by the lines
  $L_i=\langle e_i\rangle$ with $\{e_1,e_2,e_3\}$ the usual basis for
  $V$, and $c_0$ the chamber $\langle e_1\rangle\subset\langle
  e_1,e_2\rangle$ -- see Figure \ref{fig:apartment}. 
Then $P$ acts on $A_0$. Moreover,
  the isometry $\Delta_W\rightarrow A_0$, $g\mapsto gc_0$ is equivariant
  with respect to the $(W,S)$-action $g\stackrel{g_0}{\mapsto}g_0g$ on
  the Coxeter complex $\Delta_W$ and the $P$-action on the apartment
  $A_0$ (thus, the $(W,S)$-action on $\Delta_W$ is the same as the
  $P$-action on $A_0$). 
\end{description}
These three allow us to reconstruct the chambers, adjacency and
$\gS_3$-metric of $\Delta$ inside $G$:

\begin{figure}
  \centering
\begin{pspicture}(0,0)(3,3.5)
\rput(0,0.25){
\rput(1.5,1.5){\BoxedEPSF{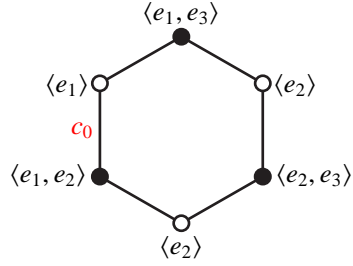 scaled 500}}
\rput(1.5,3.05){${\langle e_1, e_3\rangle}$}
\rput(3.25,0.9){${\langle e_2, e_3\rangle}$}
\rput(-0.25,0.9){${\langle e_1, e_2\rangle}$}
\rput(0,2.1){${\langle e_1\rangle}$}
\rput(1.5,-0.05){${\langle e_2\rangle}$}
\rput(3,2.1){${\langle e_2\rangle}$}
\rput(0.2,1.5){$\red{c_0}$}
}
\end{pspicture}
  \caption{Apartment $A_0$}
  \label{fig:apartment}
\end{figure}

\paragraph{Reconstructing the chambers of $\Delta$ in $G$.} 
For $a\in G$ we have $ac_0=c_0$ with 
$c_0=\langle e_1\rangle\subset\langle e_1,e_2\rangle$,
exactly when 
$$
a\in B:=
  \left\{\left(
\begin{array}{ccc}
\blob&\blob&\blob\\
0&\blob&\blob\\
0&0&\blob
\end{array}
\right)
\in G
\right\},
$$
the subgroup of upper triangular matrices. It is easy to show that
(G2) is equivalent to 
(G2a): the $G$-action on $\Delta$ is transitive on the
chambers, and (G2b): for any $g\in (W,S)$ the action of the subgroup
$B$ is transitive on the chambers $c$ such that $\delta(c_0,c)=g$.   

Combining (G2a) with the fact that the chamber $c_0$ has stabilizer
$B$, we get a 1-1 correspondence between the chambers of $\Delta$ and
the left cosets $G/B$:
$$
\xymatrix{
\text{chambers }ac_0\in\Delta\ar@{<->}[r]^-{\text{1-1}}
&  \text{cosets }aB\in G/B.
}
$$

\paragraph{Reconstructing the $i$-adjacency.} 
Let $c_1,c_2\in\Delta$ be $1$-adjacent chambers: $c_1=V_1\subset V_2$ and
$c_2=V_1'\subset V_2$, and let $c_i=a_ic_0$ with the $a_i\in
G$. Then $a_1^{-1}a_2$ stabilizes the subspace $\langle
e_1,e_2\rangle$, hence 
\begin{equation}
  \label{eq:11}
a_1^{-1}a_2\in 
  \left\{\left(
\begin{array}{ccc}
\blob&\blob&\blob\\
\blob&\blob&\blob\\
0&0&\blob
\end{array}
\right)
\in G
\right\}.  
\end{equation}
The reader can show that for $s_1$ the permutation matrix in (\ref{eq:12}),
the subgroup of matrices in (\ref{eq:11}) is the
disjoint union $B\langle s_1\rangle B:=B \cup Bs_1B$,
where $BaB=\{bab'\,|\,b,b'\in B\}$ is a double coset.
Thus, if we are to replace the chambers
$c_1,c_2$ by the cosets $a_1 B,a_2 B$, then we need to replace
$c_1\sim_1 c_2$ by  $a_1^{-1}a_2\in B\langle s_1\rangle
B$. Similarly 
$$
c_1\sim_2 c_2
\text{ exactly when the }
c_i=a_ic_0\text{ with }
a_1^{-1}a_2\in 
  \left\{\left(
\begin{array}{ccc}
\blob&\blob&\blob\\
0&\blob&\blob\\
0&\blob&\blob
\end{array}
\right)
\in G
\right\}
=B\langle s_2\rangle B.  
$$

\paragraph{Reconstructing the $\gS_3$-metric $\delta$.} 
Let $c_1,c_2\in\Delta$ be chambers with $c_i=a_ic_0$. Suppose that
$\delta(c_1,c_2)=g\in (W,S)$. As $G$ is acting by isometries
(G1), we have $\delta(c_0,a_1^{-1}a_2c_0)=g$. In the Coxeter complex
$\Delta_W$ we have by (\ref{eq:7}) 
that $\delta_W(1,g)=g$, so that by (G3),
$\delta(c_0,gc_0)=g$ also. Thus by (G2b) there is a $b\in B$ with
$(bc_0,bgc_0)=(c_0,a_1^{-1}a_2c_0)$, so in particular,
$bgc_0=a_1^{-1}a_2c_0$. As the elements of $G$ 
sending $c_0$ to $bgc_0$ are precisely the coset $bgB$, we get
$a_1^{-1}a_2\in bgB\subset BgB$. 

Conversely, if $a_1^{-1}a_2\in BgB$ then
$$
\delta(c_1,c_2)
=\delta(a_1c_0,a_2c_0)
=\delta(c_0,a_1^{-1}a_2c_0)
=\delta(c_0,bgb'c_0)
=\delta(c_0,bgc_0),
$$
for some $b\in B$, and so
$$
\delta(c_0,bgc_0)
=\delta(bc_0,bgc_0)
=\delta(c_0,gc_0)
=\delta_W(1,g)
=g,
$$
(the first as $B$ stabilizes $c_0$, the second by (G1) and the third by
(G3)). We conclude that 
$$
\delta(c_1,c_2)=g\in (W,S)\text{ if and only if }a_1^{-1}a_2\in BgB.
$$

Summarizing, let the left cosets $G/B$ be a
chamber system over $I=\{1,2\}$ with adjacency defined by
$a_1B\sim_ia_2B$ iff $a_1^{-1}a_2\in B\langle s_i\rangle B$ and
$\gS_3$-metric $\delta(a_1B,a_2B)=g$ iff $a_1^{-1}a_2\in BgB$. Then 
$G/B$ is a building of type 
\begin{pspicture}(0,0)(1.4,0.3)
\rput(0.25,0.1){
\psline(0,0)(1,0)
\pscircle[fillstyle=solid,fillcolor=white](0,0){0.125}
\pscircle[fillstyle=solid,fillcolor=white](1,0){0.125}
}
\end{pspicture},
isomorphic to the flag complex of \S \ref{lecture1}.

We leave it to the reader to show that the assumptions (G1)-(G3) hold (\emph{hint\/}: for (G2) with
$\delta(c_1,c_2)=\delta(c_1',c_2')$, situate $c_1,c_2$ in a hexagon as
in \S \ref{lecture1} and $c_1',c_2'$ similarly. Then use the fact that
$GL(V)$ acts transitively on ordered bases of $V$).

We are feeling our way towards a class of groups in which we
can mimic this reconstruction of the flag complex. 
It turns out to be
convenient to formulate the class abstractly first, and then
bring in the natural examples later. 

A \emph{Tits system\/} or 
\emph{$BN$-pair\/} for a group $G$ is a pair of subgroups $B$ and $N$
of $G$
satisfying the following axioms:
\begin{description}
\item[(BN0).] $B$ and $N$ generate $G$;
\item[(BN1).] the subgroup $T=B\cap N$ is normal in $N$, and the quotient
  $N/T$ is a Coxeter system $(W,S)$ for some $S=\{s_i\}_{i\in I}$;
\item[(BN2).] for every $g\in W$ and $s\in S$ the product of double
  cosets\footnote{A $g\in W$ is not
  an element of $G$ but a coset $\ov{g}T$ for some representative in
  $\ov{g}\in N$ for $g$. As $T\subset B$, if
  $\ov{g}_1T=\ov{g}_2T$ then $B\ov{g}_1B=B\ov{g}_2B$, so we can
  unambiguously write
  $BgB$ to mean $B\ov{g}B$.}
$BsB\cdot BgB\subset BgB\,\bigcup\, BsgB$;
\item[(BN3).] for every $s\in S$ we have $sBs\not= B$.
\end{description}
The group $W$ is called the \emph{Weyl group\/} of $G$, and is in general
not finite.

\begin{example}
\label{example:generallinear}
$G=GL_n(k)$; $B=$ the upper triangular matrices in $G$; 
$N=$ the monomial matrices in $G$ (those having exactly one non-zero entry in
each row and column), 
$$
T=\{\text{diag}(t_1,\ldots,t_n)\,|\,t_1\ldots t_n\not=0\},
$$
and $W=$ the permutation matrices with
$$
\begin{pspicture}(0,0)(\textwidth,3.75)
\rput(-0.75,-0.2){
\rput(7,2){$s_i=
\left(
  \begin{array}{cccccccc}
    1&&&&&&&\\
      &\vrule width 4mm height 0 mm depth 0mm&&&&&\\
      &&1&&&&&\\
     &&&0&1&&&\\
     &&&1&0&&&\\
     &&&&&1&&\\
     &&&&&&\vrule width 4mm height 0 mm depth 0mm&\\
     &&&&&&&1
  \end{array}
\right)$}
\rput[c](6.9,1.6){\psframe(0,0)(0.85,0.85)}
\rput(-0.6,-0.1){\psline[linewidth=1pt,linestyle=dotted](6.1,3.4)(6.9,2.9)}
\rput(0.1,0.1){\rput(2.2,-2.25){\psline[linewidth=1pt,linestyle=dotted](6.1,3.4)(6.9,2.9)}}
}
\end{pspicture}
$$
for $i\in\{1,\ldots,n-1\}$, where the number of $1$'s on
the diagonal before the $2\times 2$ block is $i-1$. 
Let $e_i$ be the $n$-column vector $(0,\ldots,1,\ldots,0)^T$ with the $1$ in the $i$-th position
and $L_i=\{te_i\,|\,t\in k\}$. Then $N$ permutes the set of lines $\{L_1,\ldots,L_n\}$ and
$W$ is isomorphic to the symmetric group on this set (hence $\cong\gS_n$).
This example is misleadingly special
in that the extension $1\rightarrow T\rightarrow N\rightarrow
W\rightarrow 1$ splits, so that the Weyl group $W$ can be realised,
via the permutation matrices, as
a subgroup of $G$. In general this doesn't happen.
\end{example}

\begin{theorem}
\label{theorem:generalized.flag.varieties}
Let $G$ be a group with a $BN$-pair and 
let $\Delta$ be a chamber system over $I$ with
chambers the cosets $G/B$ and adjacency defined by
$a_1B\sim_i a_2B$ iff $a_1^{-1}a_2\in B\langle s_i\rangle B$. Define a
$W$-metric by $\delta(a_1B,a_2B)=g\in W$ iff $a_1^{-1}a_2\in
BgB$. Then $(\Delta,\delta)$ is a thick building of type $(W,S)$.
\end{theorem}

\begin{example}
\label{example:symplectic}
$G=$ the symplectic group $\sp_{2n}(k)=\{g\in GL_{2n}(k)\,|\,g^TJg=J\}$ where 
$$
J=
\left(
  \begin{array}{cc}
    0&I_n\\
    -I_n&0
  \end{array}
\right),
$$
with $I_n$ the $n\times n$ identity matrix; $B=$ the upper triangular
matrices in $\sp_{2n}(k)$; $N=$ the monomial matrices in $\sp_{2n}(k)$, and
$$
T=\{\text{diag}(t_1,\ldots,t_n,t_1^{-1},\ldots,t_n^{-1})\,|\,t_i\not=0\}.
$$
Let $\{e_1,\ldots,e_n,\ov{e}_1,\ldots,\ov{e}_n\}$ be $2n$-column
vectors $(0,\ldots,1,\ldots,0)^T$ with the $1$ in the $i$-th
position for $e_i$ and the $(i+n)$-th position for $\ov{e}_i$.
Let $L_i=\{te_i\,|\,t\in k\}$ and
$\ov{L}_i=\{t\ov{e}_i\,|\,t\in k\}$,
writing $\ov{\ov{L}}=L$. Then $N$ permutes 
the set $\{L_1,\ldots,L_n,\ov{L}_1,\ldots,\ov{L}_n\}$
and $W$ is isomorphic to the ``signed'' permutations 
$\gS_{\pm n}=\{\pi\in\gS_{2n}\,|\,\pi(\ov{L}_i)=\ov{\pi(L_i)}\}$.

This can be reformulated geometrically as follows.
Let $V$ be a $2n$-dimensional space over $k$ and $(u,v)$ a
symplectic form on $V$ -- a non-degenerate 
alternating
bilinear form\footnote{Alternating means $(u,u)=0$ for
  all $u$, and non-degenerate that $V^\perp=\{0\}$.}. 
Let $O(V)$ be those linear maps preserving the form, i.e. $O(V)=\{g\in
GL(V)\,|\,(g(u),g(v))=(u,v)\text{ for all }u,v\in V\}$. 
The form can be defined on a basis
$\{e_1,\ldots,e_n,\ov{e}_1,\ldots,\ov{e}_n\}$ by
$$
(e_i,e_j)=0=(\ov{e}_i,\ov{e}_j)\text{ and }(e_i,\ov{e}_j)=\delta_{ij}=-(\ov{e}_j,e_i),
$$
so that $O(V)\cong \sp_{2n}(k)$.
Call a subspace $U\subset V$ totally isotropic if $(u,v)=0$
for all $u,v\in U$. 
It turns out that the maximal 
totally isotropic subspaces are $n$-dimensional. 
A (maximal) flag in $V$ is a sequence of totally
isotropic subspaces $V_1\subset\cdots\subset V_n$ with $\dim V_i=i$.
Let $\Delta$ be the chamber system with
chambers these flags and adjacencies over $I=\{1,\ldots,n\}$ as in the
flag complex of Example \ref{example:flagcomplexes}:
$(V_1\subset\cdots\subset V_{n})
\sim_i (V'_1\subset\cdots\subset V'_{n})$
when $V_j=V'_j$ for $j\not= i$. Let $c_0$ be the chamber
$$
\langle e_1\rangle
\subset 
\langle e_1,e_2\rangle
\subset\cdots\subset
\langle e_1,e_2,\ldots,e_n\rangle
$$
and $A_0$ the set of images of $c_0$ under the signed permutations
$\gS_{\pm n}=\{\pi\in\gS_{2n}\,|\,\pi(\ov{e}_i)=\ov{\pi(e_i)}\}$
(writing $\ov{\ov{e}}=e$). Finally, let $\{A_\aa\}$ be the set of
images of $A_0$ under $\sp_{2n}(k)$. Then this set of
apartments $\Delta$ gives a building isomorphic to the spherical building
of $\sp_{2n}(k)$ arising from Theorem \ref{theorem:generalized.flag.varieties} and 
Example \ref{example:symplectic}.
\end{example}

\begin{figure}
  \centering
\begin{pspicture}(0,0)(\textwidth,10)
\rput(4.65,5.25){
\rput(0,0){\BoxedEPSF{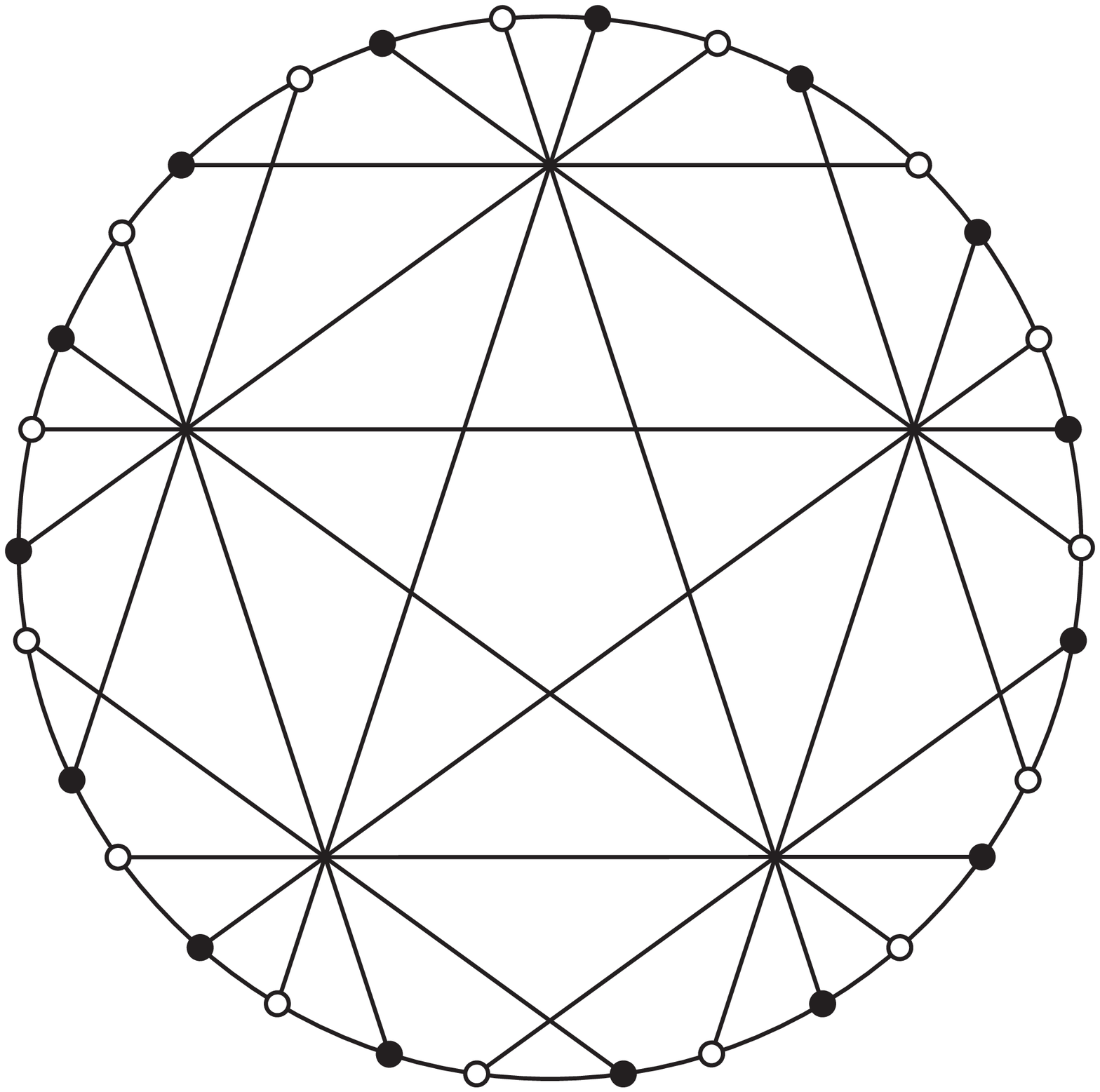 scaled 500}}
}
\rput(10.8,1.5){
\rput(0,0){\BoxedEPSF{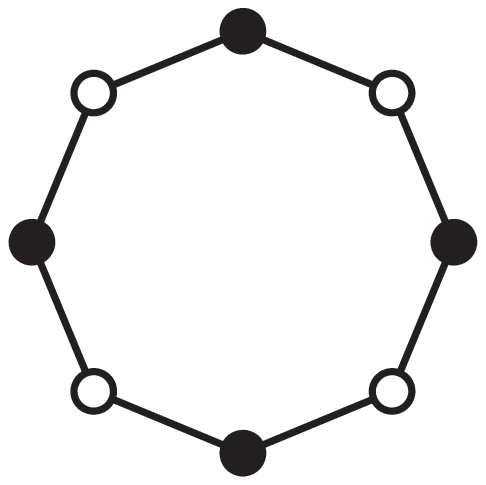 scaled 500}}
\rput(0,0){$A_0$}
\rput(1,1){$L_1$}\rput(-1,1){$L_2$}
\rput(-1,-1){$\ov{L}_1$}\rput(1,-1){$\ov{L}_2$}
\rput(0,1.35){$L_1+L_2$}\rput(0,-1.35){$\ov{L}_1+\ov{L}_2$}
\rput(1.75,0){$L_1+\ov{L}_2$}\rput(-1.75,0){$\ov{L}_1+L_2$}
}
\end{pspicture}
\caption{The spherical building of the symplectic group $\sp_4(\F_2)$
and apartment $A_0$.}
  \label{fig:flagcomplex2}
\end{figure}

We finish where we started by drawing a
picture. Let $V$ be four dimensional over the field of order $2$ and
equipped with symplectic form $(u,v)$. Let $\Delta$ be the graph with
vertices the proper non-trivial totally isotropic subspaces of $V$,
with an edge connecting the (white) one dimensional vertex $V_i$ to the
(black) two dimensional vertex
$V_j$ whenever $V_i$ is a subspace of $V_j$. Any one dimensional
subspace (of which there are 15) is totally isotropic, and is
contained in 3 two dimensional totally isotropic subspaces, each
of which in turn contains 3 one dimensional subspaces. There are
thus 15 two dimensional vertices.  The local pictures/apartments are
octagons (or barycentrically subdivided diamonds). The apartment
$A_0$ above has white vertices
$L_1,L_2,\ov{L}_1,\ov{L}_2$, using the notation of Example
\ref{example:symplectic}, and black vertices $L_1+L_2,L_1+\ov{L}_2,
\ov{L}_1 +L_2$ and $\ov{L}_1+\ov{L}_2$. 
See Figure \ref{fig:flagcomplex2}.

\begin{remark}
\label{remark:BNpair}
Examples \ref{example:generallinear} and \ref{example:symplectic} are of 
classical groups of matrices. This can be generalized.
Let $k=\ov{k}$ be 
algebraically closed and $G$ a connected algebraic group defined
over $k$. Suppose also
that $G$ is reductive, i.e. that its unipotent 
radical is trivial. Let $B$
be a Borel subgroup (a maximal closed connected soluble
subgroup) and $T\subset B$ a maximal torus --
a subgroup isomorphic to $(k^\times)^m$ for some $m$. Finally, let
$W=N/T$ be the Weyl group of $G$, where $N$ is the normalizer
in $G$ of $T$. This is isomorphic to a \emph{finite\/} Coxeter group $(W,S)$
with $S=\{s_i\}_{i\in I}$. 
The result
is a $BN$-pair for $G$. For a general non-algebraically closed $k$ 
a $BN$-pair can still be
extracted from $G$, but one has to tread more carefully.
\end{remark}

\section*{Notes and References}

As mentioned in the Introduction, most of what we have said has its
origins in the work of Tits, and we start by listing his (many)
original contributions. Coxeter groups as a notion first appeared in
his 1961 mimeographed notes, \emph{Groupes et g\'{e}om\'{e}tries de
  Coxeter\/}. These were reproduced in \cite[pages 740--754]{Wolf}. 
The name is a homage to \cite{Coxeter35}. The Bourbaki volume
\cite{Bourbaki02} dealing with Coxeter groups was produced after
``numerous conversations'' with Tits. Buildings as simplicial
complexes go back to the very beginnings of the subject, but the first
complete account can be found in \cite{Tits74}. Buildings as chamber
systems with a $W$-metric have their origins in \cite{Tits81}. The
earliest reference to $BN$-pairs that we could find in Tits's work is in
\cite{Tits62}; they start to prove an essential tool in \cite{Tits64}.

\paragraph{Section \ref{lecture1}.}
This is mostly folklore. The reader is to be minded of
projective geometry as $\Delta$ is the incidence graph of the standard
projective plane over $k$. The ad-hoc argument (essentially the
Jordan-H\"{o}lder Theorem) for associating the permutation
$(1,3)$ to the pair of chambers is from \cite[\S 4.3]{Abramenko_Brown08}.

\paragraph{Section \ref{lecture2}.}
Standard references on reflection groups and Coxeter groups are \cite{Bourbaki02}
(still the only place you can find some things), \cite{Humphreys90}
and \cite{Kane01}. The definition of reflection in (\ref{eq:1}) is
from \cite[V.2.2]{Bourbaki02}. That $\HH$ consists of all the
reflecting hyperplanes of $W$ 
is \cite[Proposition 1.14]{Humphreys90}. 
The general theory of finite reflection groups, including their
classification, can be found in Chapters 1 and 2 of \cite{Humphreys90}.
Example \ref{example:affine}, although fairly standard, is taken from
\cite[\S 2.2.2]{Abramenko_Brown08}. The general theory of affine groups
is in \cite[Chapter 4]{Humphreys90}.
For the hyperboloid or Minkowski model of hyperbolic space, hyperbolic lines, etc, see
\cite[Chapter 3]{Ratcliffe06}. The standard reference on hyperbolic
reflection groups is \cite{Vinberg85}.
The treatment of chambers, panels and adjacency
is taken from \cite[\S 1.1.4]{Abramenko_Brown08}. That $W$ acts
regularly on the chambers is \cite[Theorem 1.12]{Humphreys90}. Fact 1
is \cite[Theorem 1.5]{Humphreys90} and Fact 2 is
\cite[Theorem 1.9]{Humphreys90}. For the general theory of Coxeter
groups see
\cite[Chapter 5]{Humphreys90}. The representation 
$(W,S)\rightarrow GL(V)$ described in Remark \ref{remark:titsrep} is called
the geometric or reflectional or Tits representation, and is one of
the crucial results of \cite{Wolf}. See
\cite[\S 5.3]{Humphreys90} for its definition; faithfulness is
\cite[Corollary 5.4]{Humphreys90} or
\cite[Theorem 2.59]{Abramenko_Brown08} (where it is also shown that
the image in $GL(V)$ of $(W,S)$ is
discrete).

\paragraph{Section \ref{lecture3}.}
Apart from the aside, this section is based mainly on Chapters 1-2 of \cite{Ronan09}; the
initial chamber system notions and Example \ref{example:flagcomplexes}
are directly from \cite[\S 1.1]{Ronan09}. Chapter 2 of this book is
entirely devoted to Coxeter complexes. 
A thorough exploration of the general connections between chambers systems and
simplicial complexes is given in \cite[Appendix A]{Abramenko_Brown08}.
The building specific set-up is in \cite[\S 5.6]{Abramenko_Brown08}. 
The construction of the simplicial complex $X_\Delta$ as the
nerve of the covering by rank $|I|-1$ residues is
\cite[Exercise 5.98]{Abramenko_Brown08}.
The statement about the intersection of residues being a residue is
\cite[Exercise 5.32]{Abramenko_Brown08}. The edge coloured graph way
of viewing chamber systems is a point of view adopted in \cite{Weiss03}.

\paragraph{Section \ref{lecture4}.}
This section is based on Chapter 3 of \cite{Ronan09}
from which the definition of building is taken. That the Coxeter
complexes comprise the thin buildings is from \cite[\S 3.2]{Ronan09}. 
The alternative definition of the permutation associated to a
pair of chambers of a flag complex in
Example \ref{example:Anbuilding} is taken from \cite[Example 7.4]{Weiss03}.
The infinite 3-valent tree of Example \ref{example:affinebuilding} is
an example of a building that does not have much structure as a
combinatorial object. Nevertheless it can be constructed in an interesting
way from a vector space over a field with a discrete valuation (and as
such is an important special case of the Bruhat-Tits theory
\cite{Bruhat_Tits72}) in the following way.
Let $K$ be a non-archimedean
local field with residue field $k$ and valuation ring $A$ (for example
$K$ is the $p$-adics $\Q_p$ with $k=\Z/p\Z$ and $A$ the $p$-adic integers). If $V$ is a
$2$-dimensional vector space over $K$, then a lattice $L\subset V$ is a
free $A$-module of rank $2$. Consider the equivalence classes $\Lambda$ of
lattices under the relation $L\sim Lx$ for $x\in K^\times$, and let
$\Delta$ be the graph with vertices these classes and an edge joining
$\Lambda,\Lambda'$ iff there are $L\in\Lambda,L'\in\Lambda'$ with
$L'\subset L$ and $L/L'\cong k$. Then $\Delta$ is a tree, and Example
\ref{example:affinebuilding} is the case where $k$ has two elements
($K=\Q_2$ for example). See \cite[II.1.1]{Serre03} for
details. 
In general there is a construction that extracts a $BN$-pair, and an affine building,
from an algebraic group defined over such a $K$, and
Example \ref{example:affinebuilding} is such an affine building for
$\sl_2\Q_2$.
For affine buildings in general see \cite{Weiss09}.
The fact that the affine building for
$\sl_2\Q_p$ is a tree was used by Serre to reprove a theorem of Ihara that
a torsion free lattice in $\sl_2\Q_p$ is a free group.
A theorem of Walter Feit and Graham Higman \cite{Feit_Higman64} has
consequence that a finite thick building has type $(W,S)$ a finite
reflection group where each irreducible component of $W$ is of type
$A_n,B_n/C_n, D_n,E_6,E_7,E_8,F_4,G_2$ or $I_2(8)$
(see \cite[Theorem 6.94]{Abramenko_Brown08}; 
see \cite[Chapter 2]{Humphreys90} for a
description of these types of finite reflection group). 
Hence there can be no finite thick buildings of type the
symmetry group of the dodecahedron, for which $(W,S)$ has type
$H_3$. That there are no \emph{infinite\/} thick buildings of type
$H_3$ is shown in \cite{Tits77}.
Theorem \ref{theorem:isometries} is
\cite[Theorem 3.6]{Ronan09} and Theorem
\ref{theorem:second_building_definition} is \cite[Theorem
3.11]{Ronan09}. Prior to \cite{Tits81} axioms (B1$^\prime$) and (B2$^\prime$) of
Theorem \ref{theorem:second_building_definition} provided the standard
definition of building.

\paragraph{Section \ref{lecture5}.}
This section is based on \cite[Chapter 5]{Ronan09}. Properties
(G1)-(G3) are the specialization to $GL_3$ of a strongly transitive
group action \cite[\S 5.1]{Ronan09}. The argument that reconstructs
the $W$-metric is taken from the proof of \cite[Theorem 5.2]{Ronan09}. 
The axioms for a $BN$-pair are from
\cite[\S5.1]{Ronan09}. 
A proof that Example \ref{example:generallinear} is a $BN$-pair using
nothing but row and column operations can be found in \cite[\S 6.5]{Abramenko_Brown08}.
Theorem
\ref{theorem:generalized.flag.varieties} is \cite[Theorem 5.3]{Ronan09}. 
The flag complex of a symplectic space is from
\cite[Chapter 1]{Ronan09}.
Figure \ref{fig:flagcomplex2} has several names: in graph theory circles it is called
Tutte's eight-cage, and is the unique smallest cubic graph with girth
$8$ (where these minimal $8$-circuits are, of course, the
apartments). It is a pleasantly mindless exercise to label the vertices of the
Figure with the totally isotropic subspaces (\emph{hint:\/} start with
the $8$-circuit at the top as the apartment $A_0$). 
There is also a very simple construction that goes back to Sylvester
(1844) -- this (and much else) is engagingly described in
\cite{Coxeter58}. There are $30$ odd permutations of order $2$ in
$\gS_6$: $15$ transpositions -- like $(1,2)$ -- and $15$ products of three
disjoint transpositions, like $(1,2)(3,4)(5,6)$. Let these be the
vertices of the eight-cage, and join a vertex $\ss$ in one of these
two groups to the three $\tau_1,\tau_2,\tau_3$ in the other group for which
$\ss=\tau_1\tau_2\tau_3$. 
That the $B$
(Borel subgroup) and $N$ (normalizer of a maximal torus) extracted
from a reductive group $G$ in Remark \ref{remark:BNpair} are a
$BN$-pair for $G$ is shown in
\cite[\S 29.1]{Humphreys75}. 

\paragraph{Further reading.}
Surely the shortest introduction to buildings is
\cite{Brown02}; \cite{Brown91},
\cite{Rousseau09} and \cite{Tits75} are slightly longer.
The book \cite {Abramenko_Brown08} is a greatly expanded version of 
\cite{Brown89}, while \cite{Ronan09} is an updated version of the 1988 original.
A nice introduction to spherical buildings, including an account of
Tits's classification \cite{Tits74} of the thick spherical buildings of
type $(W,S)$ for $|S|\geq 3$, is
\cite{Weiss03}; the sequel \cite{Weiss09} treats affine buildings.




\bibliography{version3}{}
\bibliographystyle{plain}



%
%
%
%

\end{document}